\documentclass{article}

\usepackage[utf8]{inputenc}
\usepackage{amsmath}
\usepackage{amsfonts}
\usepackage{amssymb}
\usepackage{stmaryrd}

\usepackage{pgfplots}
\usepackage{tikz}
\pgfplotsset{height=8cm, width=15cm, compat=1.17}

\usepackage[ngerman]{babel}

\usepackage{url}

\usepackage[a4paper,top=2.5cm,bottom=2cm,left=2.5cm,right=2.5cm,marginparwidth=1.75cm]{geometry}

\addto\captionsngerman{}

\begin{document}

\begin{titlepage}
    \begin{center} 
        \vspace*{\fill}
    
        \Huge \textbf{Über ein kombinatorisches Problem mit $n$ Plätzen und $n$ Personen
        }
        \vspace{4cm}

        \Large Simon Wundling\\
        \Large swundling@gmx.de \\
        \Large Gymnasium der Benediktiner, Meschede\\
        
        \vspace{2cm}

        \Large 31.03.2023
    
        \vspace{2.5cm}

        \begin{abstract}
            \noindent Möchte man $n \in \mathbb{N}$ Plätze nacheinander mit $n$ Personen und der Regel, dass jede Person einen der Plätze mit maximalen Abstand zu einem besetzten Platz auswählt, besetzen, dann kann man sich fragen, wie viele Möglichkeiten es dafür gibt. In dieser Arbeit soll anhand anfänglich genannter Ideen eine Formel für die Anzahl dieser Möglichkeiten gefunden werden. Zudem wird diese Formel nach unten und oben abgeschätzt. Schlussendlich werden Formeln für die OEIS-Folgen A166079, A095236, A095240 und A095912 und eine Erweiterung des anfänglichen Problems hergeleitet.
        \end{abstract}

        \vfill
    \end{center}
\end{titlepage}

\newpage

\tableofcontents

\begin{center}
    
\end{center}

\vspace{0.5cm}

\begin{Large}
    \noindent \textbf{Anderes} \\
\end{Large}

\begin{itemize}
    \item Definition Teilreihe: Reihe von benachbarten leeren Plätzen (Die Anzahl der leeren Plätze gibt die Länge einer Teilreihe an)
    \item Im Folgenden soll mit Abstand immer der Abstand zu einem besetzten Platz gemeint sein (Bei abweichender Verwendung wird explizit genannt, was gemeint ist).
    \item X als Symbol für einen besetzten Platz
    \item O als Symbol für einen leeren Platz
\end{itemize}

\newpage

\section{Einführung}

\noindent Sei eine Reihe von $n \in \mathbb{N}$ Plätzen gegeben (alternativ können auch $n$ Pissoirs wie in [1] betrachtet werden). Benachbarte Plätze haben einen Abstand von $1$. Nun sollen $n$ Personen nacheinander einen der $n$ Plätze besetzen und sich dabei an folgende Regel halten:\\

\noindent \textbf{Höflichkeitsregel:} Wer dazu kommt, wählt einen der Plätze so aus, dass der Abstand zu einem besetzten Platz maximal wird. \\

\noindent Bei $n=8$ und der ersten Person auf dem Platz links außen folgt exemplarisch:

\begin{enumerate}
    \item XOOOOOOO
    \item XOOOOOOX
    \item XOOXOOOX oder XOOOXOOX
    \item XOOXOXOX oder XOXOXOOX
    \item Besetzen der verbliebenen leeren Plätze (alle haben den Abstand $1$)
\end{enumerate}

\noindent Nun kann man die Frage stellen, wie viele Möglichkeiten $a_n$ es gibt, dass $n$ Personen nacheinander jeweils einen der $n$ Plätze besetzen und sich dabei an die Höflichkeitsregel halten. Bei sehr kleinen $n$ ist die Berechnung im Kopf oder per Hand noch möglich. Es ist $a_1 = 1$, da es nur eine Möglichkeit (X) gibt, und $a_2=2$, da folgende zwei Möglichkeiten existieren:
\begin{center}
    1. XO und dann XX \ \ \ \ \ \ \ \ \ \ \ 2. OX und dann XX
\end{center}

\noindent Wird aber z.B. $n=10$ betrachtet, dann ist $a_{10}$ sehr schwierig per Hand oder im Kopf zu bestimmen. In dieser Arbeit soll es deswegen darum gehen, eine Formel für $a_n$ zu finden und diese zu beweisen. Dafür muss also untersucht werden, welche Aspekte für die Berechnung von $a_n$ relevant sind.  Zunächst gibt es, sofern mindestens ein Platz besetzt ist, zwei charakteristische Arten, auf die die nächste Person einen Platz besetzen kann:

\begin{align*}
    1. \ \ X\underbrace{\text{O...O}}_{\text{Länge $L$}}
    \ \ \ \ \ \ \ \ \ \ \
    2. \ \ X\underbrace{\text{O...O}}_{\text{Länge $L$}}X
\end{align*}

\noindent Beim ersten Fall besetzt dann die nächste Person den äußeren Platz der Teilreihe mit Abstand $L$. Dagegen wird beim zweiten Fall der Platz bzw. einer der Plätze in der Mitte der Teilreihe besetzt (Abstand $\left \lfloor \frac{L-1}{2} \right \rfloor + 1$). Bezüglich des zweiten Falls muss beim Besetzen zwischen geradem und ungeradem $L$ unterschieden werden. Ist nämlich $L$ ungerade, dann gibt es für die nächste Person nur eine Möglichkeit, einen Platz in der Mitte der Teilreihe zu besetzen. Bei geradem $L$ gibt es hingegen zwei Plätze, die als Mitte der Teilreihe in Frage kommen, womit die entsprechende Person zwischen zwei Möglichkeiten auswählen kann. Der nicht ausgewählte Platz hat dann einen Abstand von $1$ und spielt somit für die Besetzung von Plätzen mit Abstand $\geq 2$ keine Rolle. Des Weiteren gibt es dabei den Spezialfall $L=2$. Dort fällt der Aspekt der zwei Möglichkeiten weg, da dieser schon durch die Anzahl der Personen, die sich auf einen Platz mit Abstand 1 setzen, dargestellt wird.\\

\noindent Sei $T_j \ (j \in \mathbb{N})$ die $j$-te Teilreihe von links.  Innerhalb von $T_j$ gibt es nun einen Platz bzw. zwei Plätze mit dem Abstand $D_j$, der größer als die Abstände der anderen Plätze von $T_j$ ist. 
Die nächste Person besetzt dann einen Platz mit Abstand $H = \text{max} \lbrace D_1,D_2,D_3,... \rbrace$. 
Es entstehen eine bzw. zwei neue Teilreihen, deren jeweiliger charakteristischer Abstand $D'$ bzw. $D'$ und $D''$ kleiner als $\text{max} \lbrace D_1,D_2,D_3,... \rbrace$ ist.
Nun kann man das Prinzip mit $ \lbrace \ (\lbrace D_1,D_2,D_3,... \rbrace \setminus \lbrace H \rbrace ) \ \cup \ \lbrace D' \rbrace \ \rbrace$ bzw. $ \lbrace \ (\lbrace D_1,D_2,D_3,... \rbrace \setminus \lbrace H \rbrace ) \ \cup \ \lbrace D' , D'' \rbrace \ \rbrace$ wiederholen.
\\

\noindent Somit werden erst die Plätze mit Abstand $n-1$ ($\text{max} \lbrace D_1,D_2,D_3,... \rbrace$ mit maximalem Wert $n-1$) und dann die Plätze mit Abstand $n-2/n-3/.../2/1$ besetzt. \\

\noindent Insgesamt gibt es also zwei für die Berechnung von $a_n$ relevante Aspekte:
\begin{enumerate}
    \item Anzahl der Personen, die einen Platz mit Abstand $k \in \mathbb{N}$ besetzen
    \item Anzahl der unterschiedlichen Mengen zweier benachbarter Plätze, bei denen beide Plätze zu einem Zeitpunkt den Abstand $k \geq 2$ besitzen
\end{enumerate}

\noindent Diese beiden Aspekte werden im Folgenden mit der Einführung von $b_{p,k}$ und $d_{p,k}$ untersucht.

\section{Schemata für das Hinzufügen}

\noindent Hat man bei einer Situation mit $n' \in \mathbb{N}$ Plätzen, bei der die erste Person den Platz links außen besetzt, ermittelt, wie viele Personen einen Platz mit Abstand $k \in \mathbb{N}$ besetzen, dann stellt sich die Frage wie viele Personen es bei $n'+1/n'+2/n'+3/...$ sind. Hat man zum Beispiel bei $n'=5$ mit XOXOX eine Person, die einen Platz mit Abstand 2 besetzt, dann ist die Frage wie viele Personen bei $6/7/8/...$ Plätzen einen Platz mit Abstand 2 besetzen.\\

\noindent Dazu könnte man die im Vergleich zu $n'$ hinzugekommenen leeren Plätze in die vorhandenen Teilreihen so einfügen, dass im Nachhinein alle schon besetzten vorhandenen Plätze gemäß der Höflichkeitsregel besetzt wurden. Aus den daraus entstehenden Situationen könnte man dann die Anzahl der Personen bestimmen. Hierfür soll nun ein Schema angegeben und bewiesen werden.\\

\noindent Sei für alle $i \in \mathbb{N}$ $S_{2^i} = (s_{2^i,1};s_{2^i,2};...;s_{2^i,2^i})$ und $S_2 = (1;2)$. Weiter soll gelten: \\
\begin{equation*}
     s_{2^{i+1},k} =
    \begin{cases}
        2 \cdot s_{2^i, k } - 1
        & , \text{wenn} \ k \leq 2^i \\
        
         2 \cdot s_{2^i, k-2^i}
        & , \text{wenn} \ k > 2^i \\
        \end{cases}
\end{equation*}

\noindent Zum Beispiel ist $S_4 = (1;3;2;4)$ und
$S_8 = (1;5;3;7;2;6;4;8)$. \\

\noindent Sei nun eine Situation, bei der die erste Person den Platz links außen besetzt und $2^h \ (h \in \mathbb{N})$ Teilreihen einer bestimmten Länge $l \in \mathbb{N}$ vorhanden sind, gegeben. Somit gibt es $2^{h-1}$ besetzte Plätze mit Abstand $l+1$. Besetzt man nun Plätze, dann tut man dies immer, indem man einen Platz in der Mitte einer Teilreihe auswählt. Dass dabei die Anzahl der Plätze in der entstehenden linken Teilreihe von der Anzahl der Plätze in der entstehenden rechten Teilreihe höchstens um 1 abweicht, ist ein zentraler Aspekt.\\

\noindent Weiter weisen wir jeder Teilreihe von links nach rechts eine Zahl zu. Dabei soll die erste Teilreihe von links die Zahl $1$, die zweite Teilreihe von links die Zahl $2$, ..., die $2^h$-te Teilreihe von links die Zahl $2^h$ erhalten. Falls dabei eine Teilreihe mit einer Zahl durch Besetzen aufgespalten wird, sollen diesen zwei neuen Teilreihen zunächst keine Zahl zugewiesen werden. Erst wenn $2^{h+1}$ Teilreihen vorhanden sind, soll allen Teilreihen wieder eine Zahl zugewiesen werden. \\

\noindent \textbf{Schema 1 (mindestens eine Länge einer Teilreihe im Intervall $[l;2l)$):}
\begin{enumerate}
    \item Der erste leere Platz wird in die Teilreihe mit der Zahl $s_{2^i,1}$ eingefügt
    \item Wurde zuletzt in die Teilreihe mit der Zahl $s_{2^i,c}\ (c \in \mathbb{N} \ \text{und} \ c < 2^i)$ ein leerer Platz eingefügt, dann wird der nächste leere Platz in die Teilreihe mit der Zahl $s_{2^i,c+1}$ eingefügt
    \item Wurde zuletzt in die Teilreihe mit der Zahl $s_{2^i,2^i}$ ein leerer Platz eingefügt, dann wird der nächste leere Platz in die Teilreihe mit der Zahl $s_{2^i,1}$ eingefügt
\end{enumerate}

\noindent Dass dieses Schema tatsächlich funktioniert, kann man mit dem Prinzip der vollständigen Induktion zeigen. Zudem reicht es aus, alles für einen Durchlauf (Nach einem Durchlauf wurde in jede Teilreihe ein Platz hinzugefügt) zu zeigen. Nach einem Durchlauf beginnt das Verfahren nämlich wieder von vorn (siehe Punkt 3).\\

\noindent Zunächst kann man $2^1$ Teilreihen betrachten. Damit der besetzte Platz zwischen diesen beiden Teilreihen im Nachhinein korrekt besetzt wurde, darf die Differenz aus den Längen der Teilreihen höchstens 1 betragen. Dies erreicht man, indem man mit Schema 1 (mit $S_{2^1}$) vorgeht. Ein abwechselndes Hinzufügen von leeren Plätzen in die beiden Teilreihen führt nämlich dazu, dass jene Differenz nie größer als 1 wird. \\

\noindent Nun kann man für ein $h \in \mathbb{N}$ annehmen, dass bei $2^h$ Teilreihen Schema 1 mit $S_{2^h}$ genügt, um die leeren Plätze korrekt hinzuzufügen. Um die Induktion abzuschließen, fehlt nur noch der Schritt, dass man ausgehend von dieser Annahme zeigt, dass bei $2^{h+1}$ Teilreihen Schema 1 mit $S_{2^{h+1}}$ genügt. Zuerst stellt man folgendes fest:

\begin{itemize}
    \item Teilreihe mit Zahl $j \in \mathbb{N}$ bei $2^h$ Teilreihen ($S_{2^h}$) wird zu den Teilreihen mit den Zahlen $2j-1$ und $2j$ bei $2^{h+1}$ Teilreihen ($S_{2^{h+1}}$)
\end{itemize}

\noindent Grund dafür ist, dass man $2^{h+1}$ Teilreihen erhält, indem man $2^h$ Teilreihen jeweils in der Mitte durch das Besetzen eines Platzes spaltet. Nun sind für diese neu entstandenen besetzten Plätze nur die Differenz der Längen der jeweils erzeugten zwei Teilreihen relevant. Damit alle anderen besetzten Plätze im Nachhinein korrekt besetzt wurden, kann das Prinzip von $S_{2^h}$ verwendet werden. Wichtig ist hierbei nur, dass unter den ersten $2^h$ Teilreihen, in die ein leerer Platz hinzugefügt wird, sich nicht sowohl die Teilreihe mit der Zahl $2j-1$ als auch die Teilreihe mit der Zahl $2j$ befinden.\\

\noindent Somit reicht es aus, das Prinzip von $S_{2^h}$ zweimal folgendermaßen anzuwenden: einmal mit den Zahlen $1,3,5,...$ und einmal mit den Zahlen $2,4,6,...$ . Man erhält dann folgendes Tupel $S_{2^{h+1}}$:
\begin{equation*}
    S_{2^{h+1}} =
    (2\cdot s_{2^h,1}-1 ; 2 \cdot s_{2^h,2}-1 ;...; 2 \cdot s_{2^h,2^h}-1 
    ; 2\cdot s_{2^h,1} ; 2 \cdot s_{2^h,2} ;...; 2 \cdot s_{2^h,2^h})
\end{equation*}

\noindent Da das Schema bei $2^1$ Teilreihen gültig ist und man aus der Gültigkeit von $S_{2^h}$ die Gültigkeit von $S_{2^{h+1}}$ schließen kann, ist die Induktion fertig. Außerdem ist:

\begin{equation*}
     s_{2^{h+1},k} =
    \begin{cases}
        2 \cdot s_{2^h, k } - 1
        & , \text{wenn} \ k \leq 2^h \\
        
         2 \cdot s_{2^h, k-2^h}
        & , \text{wenn} \ k > 2^h \\
        \end{cases}
\end{equation*}

\noindent Irgendwann ist dann die Situation gegeben, dass die Längen aller Teilreihen gleich $2l$ sind. Nun kann man in die Mitte einer solchen Teilreihe einen besetzten Platz mit Abstand $l+1$ einfügen. Dafür führen wir auch ein Schema ein: \\

\noindent \textbf{Schema 2 (Längen aller $2^h$ Teilreihen gleich $2l$):}
\begin{enumerate}
    \item Der erste besetzte Platz wird in die Mitte der Teilreihe mit der Zahl $s_{2^i,1}$ eingefügt
    \item Wurde zuletzt in die Mitte der Teilreihe mit der Zahl $s_{2^i,c} \ (c \in \mathbb{N} \ \text{und} \ c < 2^i)$ ein besetzter Platz eingefügt, dann wird der nächste besetzte Platz in die Mitte der Teilreihe mit der Zahl $s_{2^i,c+1}$ eingefügt
\end{enumerate}

\noindent Da es für die Einhaltung der Höflichkeitsregel bzw. das Hinzufügen von Plätzen egal ist, ob ein leerer oder ein besetzter Platz hinzugefügt wurde (Es kann auch erst ein leerer Platz eingefügt werden, der dann später besetzt wird), sind wir nach der zuvorigen Argumentation bei Schema 1 hier fertig.

\section{Berechnung von $b_{p,k}$}

\noindent Sei $b_{p,k}$ die Anzahl der Personen, die einen Platz mit Abstand $k \in \mathbb{N}$ besetzen, wenn $p \in \mathbb{N}$ Personen nacheinander gemäß der Höflichkeitsregel einen der $p$ Plätze besetzen und die erste Person den Platz links außen besetzt.

\subsection{Berechnung von $b_{p,k}$ für $k \geq 2$}

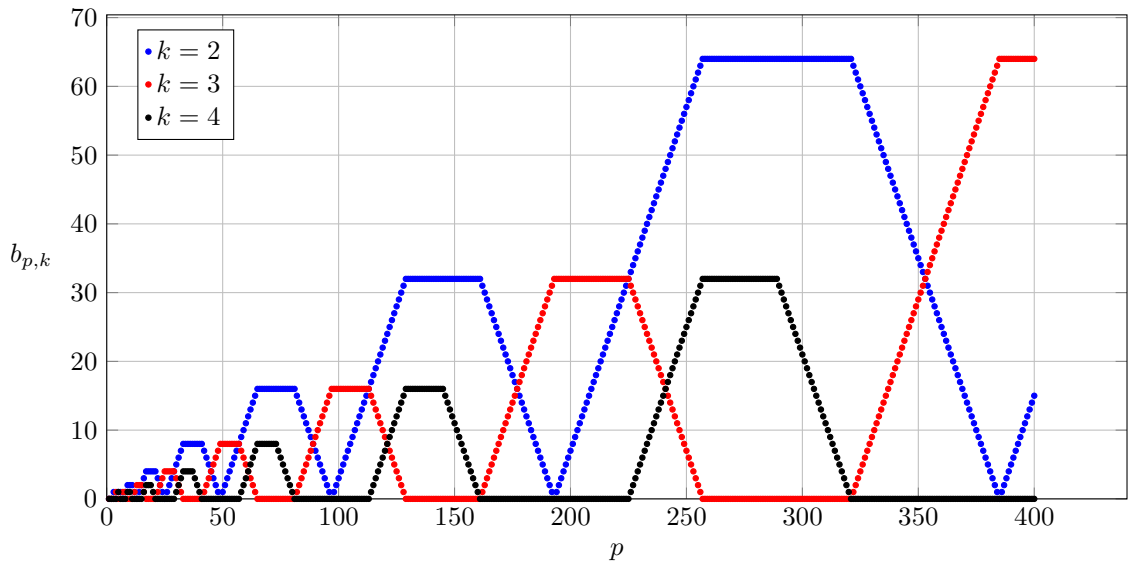
\begin{center}
    \begin{tikzpicture}
        \begin{axis} [xlabel={$p$}
        , ylabel={$b_{p,k}$}
        , ylabel style={rotate=-90}
        , xmin=0
        , ymin = 0
        , grid=major
        , legend pos=north west]
            \addplot[only marks,mark=*, mark size = 1pt, color=blue] table[col sep=semicolon]{b_n,k/b_n,2.csv};
            \addlegendentry{$k=2$}
            \addplot[only marks,mark=*, mark size = 1pt, color=red] table[col sep=semicolon]{b_n,k/b_n,3.csv};
            \addlegendentry{$k=3$}
            \addplot[only marks,mark=*, mark size = 1pt, color=black] table[col sep=semicolon]{b_n,k/b_n,4.csv};
            \addlegendentry{$k=4$}
        \end{axis}
    \end{tikzpicture}
\end{center}

\noindent \textbf{Behauptung 3.1:} Mit $m = \lfloor \log_{2}(\frac{p-1}{2k}) \rfloor$ für $p \geq 5$ gilt für alle $k \geq 2$ und alle $p$
\begin{equation*}
    b_{p,k} = 
    \begin{cases}
        0
        \ \ & \text{für} \ \ p < k+1  \\
        
        1
        \ \ & \text{für} \ \ p = k+1  \\
        
        0
        \ \ & \text{für} \ \ k+1< p < 1+2k  \\
    
        2^m 
        \ \ & \text{für} \ \ 1+2^{m}\cdot 2k \leq p \leq 1+2^{m}(2k+1) \\

        1+2^m(2k+2)-p
        \ \ & \text{für} \ \ 1+2^{m}(2k+1) < p \leq 1+2^{m}(2k+2)\\
        
        0  
        \ \ & \text{für} \ \ 1+2^{m}(2k+2)< p \leq 1+2^{m}(4k-2) \\
        
        p-1-2^{m}(4k-2)
        \ \ & \text{für} \ \ 1+2^{m}(4k-2) < p < 1+2^{m+1}2k

    \end{cases}
\end{equation*}

\noindent \textbf{Beweis:} Im Folgenden sei $k\geq 2$. Zudem soll immer die Situation, dass die erste Person den Platz links außen besetzt, betrachtet werden. Für $p \leq k$ ist $b_{p,k}=0$, da der größtmögliche Abstand (Abstand vom Platz links außen zum Platz rechts außen) höchstens $p-1$ beträgt. Man erhält:
\begin{equation}
    b_{p,k} = 0 \ \ \text{für} \ \ p \leq k
\end{equation}

\noindent Ist $p=k+1$, dann beträgt der Abstand vom Platz links außen zum Platz rechts außen $k$, womit zunächst $b_{p,k}=1$. Es bleibt eine Teilreihe der Länge $k-2$ übrig. Die dritte Person besetzt dann gemäß der Höflichkeitsregel den Platz in der Mitte dieser Teilreihe. Dieser Platz besitzt den Abstand $\left \lfloor \frac{k-2-1}{2} \right \rfloor + 1 < k$. Somit ist der anfangs erwähnte Platz der einzige Platz mit Abstand $k$:
\begin{equation}
    b_{p,k} = 1 \ \ \text{für} \ \ p = k+1
\end{equation}

\noindent Wenn $k+1<p<2k+1$, dann bleiben, nachdem die erste Person den Platz links außen und die zweite Person den Platz rechts außen besetzt hat (Abstand dieser beiden Plätze ist $\geq k+1$), noch eine Teilreihe, die höchstens eine Länge von $2k-2$ hat, übrig. Somit hat der Platz, den die dritte Person gemäß der Höflichkeitsregel in der Mitte dieser Teilreihe besetzt, höchstens einen Abstand von $\lfloor \frac{2k-2-1}{2}\rfloor +1= k-1$. Folglich gilt:
\begin{equation}
    b_{p,k} = 0 \ \ \text{für} \ \ k+1 < p < 2k+1
\end{equation}

\noindent Im Folgenden wollen wir für alle $m \in \mathbb{N}_0$ eine allgemeine Herleitung der Werte von $b_{p,k}$ im Intervall $[1+2^{m}\cdot 2k \ ;1+2^{m+1}\cdot 2k)$ durchführen. Zudem ist der Abstand des Platz links außen zum Platz rechts außen für den Wert von $b_{p,k}$ nicht mehr relevant, da dieser mindestens $1+2^0\cdot 2k -1= 2k$ beträgt. Von nun an gehen wir deshalb immer von der Situation aus, dass diese beiden Plätze bereits besetzt sind. Ferner sollen die Schemata aus Abschnitt 2 für das korrekte Hinzufügen von leeren oder besetzten Plätzen verwendet und immer wieder fortgesetzt werden. \\

\noindent Sei $p=1+2^{m}\cdot 2k$. Zunächst sind ja schon der linke und der rechte äußere Platz besetzt, womit eine Teilreihe der Länge $1+2^{m}\cdot 2k-2=2^{m}\cdot 2k-1$ entsteht. Gemäß der Höflichkeitsregel setzt sich nun die dritte Person in die Mitte dieser Teilreihe, womit zwei andere Teilreihen der Länge $\frac{(2^{m}\cdot 2k-1)-1}{2} = 2^{m-1}\cdot 2k -1$ entstehen. Nun können wir je eine dieser Teilreihen erneut in zwei Teilreihen durch das Besetzen eines Platzes in der Mitte dieser Teilreihe aufspalten und so immer wieder mit jeder Teilreihe vorgehen. Man erhält also zu bestimmten Zeitpunkten:
\begin{itemize}
    \item zwei Teilreihen der Länge $\frac{(2^{m}\cdot 2k-1)-1}{2} = 2^{m-1}\cdot 2k -1$ 
    \item vier Teilreihen der Länge $\frac{(2^{m-1}\cdot 2k-1)-1}{2}  = 2^{m-2}\cdot 2k -1$ 
    \item ...
    \item $2^{m}$ Teilreihen der Länge $\frac{(2^{1}\cdot 2k-1)-1}{2} = 2^{0}\cdot 2k -1 = 2k-1$ 
\end{itemize}

\noindent Nun kann man in die Mitte jeder Teilreihe der Länge $2k-1$ einen Platz besetzen. Diese Plätze haben dann alle einen Abstand von $\frac{(2k-1)-1}{2}+1 = k$. Da es $2^{m}$ Teilreihen der Länge $2k-1$ gibt, ist $b_{p,k}=2^m$ bei $p=1+2^{m}\cdot 2k$. Zudem gibt es nun $2^{m+1}$ Teilreihen der Länge $k-1$. \\

\noindent Nun hat jeder besetzte Platz mit Abstand $k$ links und rechts von sich eine Teilreihe der Länge $k-1$. Nun kann man für jeden diesen Platz in eine der beiden benachbarten Teilreihen einen weiteren leeren Platz hinzufügen, ohne dass sich der Abstand vom besetzten Platz verändert. Da man in dieser Weise $2^m$ leere Plätze nach dem Schema 1 (mit $S_{2^{m+1}}$) hinzufügen kann, folgt:
\begin{gather}
    b_{p,k} = 2^m \nonumber \\ 
    \text{für} \\ 
    1+2^m \cdot 2k \leq p \leq 1+2^m \cdot 2k +2^m = 1+2^m(2k+1) \nonumber
\end{gather}

\noindent Nun kann man weitere leere Plätze hinzufügen. Pro neuen hinzugefügtem Platz entsteht eine neue Teilreihe der Länge $k$. Da wir bei $p=1+2^m(2k+1)$ schon $2^m$ Teilreihen der Länge $k$ haben, haben dann pro neuer entstandener Teilreihe der Länge $k$ ein jeweils neuer Platz links und rechts von sich eine Teilreihe der Länge $k$. Der Abstand von einem solchen Platz beträgt dann $k+1$, womit pro neuem Platz dieser Art $b_{p,k}$ um $1$ sinkt. Somit ist für die Berechnung von $b_{p,k}$ der Abstand zu $1+2^m(2k+1)$ relevant, da dieser angibt, zu wie vielen Teilreihen der Länge $k-1$ ein Platz hinzugefügt wurde. Bei $p=1+2^m(2k+1)$ gibt es $2^m$ Teilreihen der Länge $k-1$. Man kann folglich $2^m$ leere Plätze mit Fortsetzung des Schemas 1 hinzufügen, bis überall Teilreihen der Länge $k$ vorhanden sind. Somit:
\begin{gather}
    b_{p,k} = 2^m - (p-(1+2^m(2k+1)) = 1+2^m(2k+2)-p \nonumber \\ 
    \text{für} \\ 
    1+2^m(2k+1) < p \leq 1+2^m(2k+1) +2^m = 1+2^m(2k+2) \nonumber
\end{gather}

\noindent Nun liegen $2^{m+1}$ Teilreihen der Länge $k$ vor. Jetzt kann man weitere leere Plätze mit Fortsetzung des Schemas 1 hinzufügen bis überall Teilreihen der Länge $2k-2$ sind. Wenn man irgendwann in der Mitte einer beliebigen Teilreihe einen Platz besetzen würde, dann hätte dieser höchstens einen Abstand von $\left \lfloor \frac{2k-2-1}{2} \right \rfloor + 1 = k-1$. Folglich ist:
\begin{gather}
    b_{p,k} = 0 \nonumber \\ 
    \text{für} \\ 
    1+2^m(2k+2) < p \leq 1+2^m(2k+2) +2^{m+1}(k-2) = 1+2^m(4k-2) \nonumber
\end{gather}

\noindent Nun haben alle Teilreihen die Länge $2k-2$. Fügt man in eine dieser Teilreihen einen neuen Platz hinzu, dann liegt eine Teilreihe der Länge $2k-1$ vor. Besetzt man in der Mitte einer solchen Teilreihe einen Platz, dann hat dieser einen Abstand von $\lfloor \frac{2k-1}{2} \rfloor +1 = k$, womit sich $b_{p,k}$ um $1$ erhöht. Nun kann man dies nach Schema 2 (mit $S_{2^{m+1}}$) auf alle $2^{m+1}$ Teilreihen der Länge $2k-2$ anwenden. Somit wird der Abstand zu $p=1+2^m(4k-2)$ relevant, da dieser angibt, zu wie vielen Teilreihen der  Länge $2k-2$ ein weiterer besetzter Platz hinzugefügt wurde. Es ist:
\begin{gather}
    b_{p,k} = p-(1+2^m(4k-2))=p-1-2^m(4k-2) \nonumber \\ 
    \text{für} \\ 
    1+2^m(4k-2) < p < 1+2^m(4k-2) +2^{m+1} = 1+2^{m+1} \cdot 2k \nonumber
\end{gather}

\noindent Mit (1) bis (3) kann man $b_{p,k}$ im Intervall $[1;2k]$ berechnen. Mit (4) bis (7) kann man $b_{p,k}$ für das Intervall $[1+2^{m}\cdot 2k \ ;1+2^{m+1}\cdot 2k)$ für alle $m \in \mathbb{N}_0$ berechnen. Da $[1;2k] \cup [1+2^{0}\cdot 2k \ ;1+2^{1}\cdot 2k) 
\cup [1+2^{1}\cdot 2k \ ;1+2^{2}\cdot 2k)
\cup [1+2^{2}\cdot 2k \ ;1+2^{3}\cdot 2k) 
\cup ... = [1; \infty)$, kann man also mit (1) bis (7) $b_{p,k}$ für alle $p$ berechnen. \\

\noindent Zudem gilt:
\begin{align*}
    p &= 1+2^{m}\cdot 2k \\
    \Leftrightarrow \frac{p-1}{2k} &= 2^m \\
    \Leftrightarrow m &= \log_{2} \left ( \frac{p-1}{2k} \right ) 
\end{align*}

\noindent Da $m$ im Intervall $[1+2^{m}\cdot 2k \ ;1+2^{m+1}\cdot 2k)$ den gleichen ganzzahligen Wert annehmen soll, muss $m = \lfloor \log_{2}(\frac{p-1}{2k}) \rfloor$. Mit $1+2^0 \cdot 4 = 5$ als kleinste Zahl, die in einem solchen Intervall liegt, sei $m = \lfloor \log_{2}(\frac{p-1}{2k}) \rfloor$ für $p \geq 5$.\\

\noindent Mit (1),(2),...,(7) folgt schließlich:\\
\begin{equation*}
    b_{p,k} = 
    \begin{cases}
        0
        \ \ & \text{für} \ \ p < k+1  \\
        
        1
        \ \ & \text{für} \ \ p = k+1  \\
        
        0
        \ \ & \text{für} \ \ k+1<p < 1+2k  \\
    
        2^m 
        \ \ & \text{für} \ \ 1+2^{m}\cdot 2k \leq p \leq 1+2^{m}(2k+1) \\

        1+2^m(2k+2)-p
        \ \ & \text{für} \ \ 1+2^{m}(2k+1) < p \leq 1+2^{m}(2k+2)\\
        
        0  
        \ \ & \text{für} \ \ 1+2^{m}(2k+2)< p \leq 1+2^{m}(4k-2) \\
        
        p-1-2^{m}(4k-2)
        \ \ & \text{für} \ \ 1+2^{m}(4k-2) < p < 1+2^{m+1}2k

    \end{cases}
\end{equation*}
\hfill $\square$
 
\subsection{Berechnung von $b_{p,1}$}

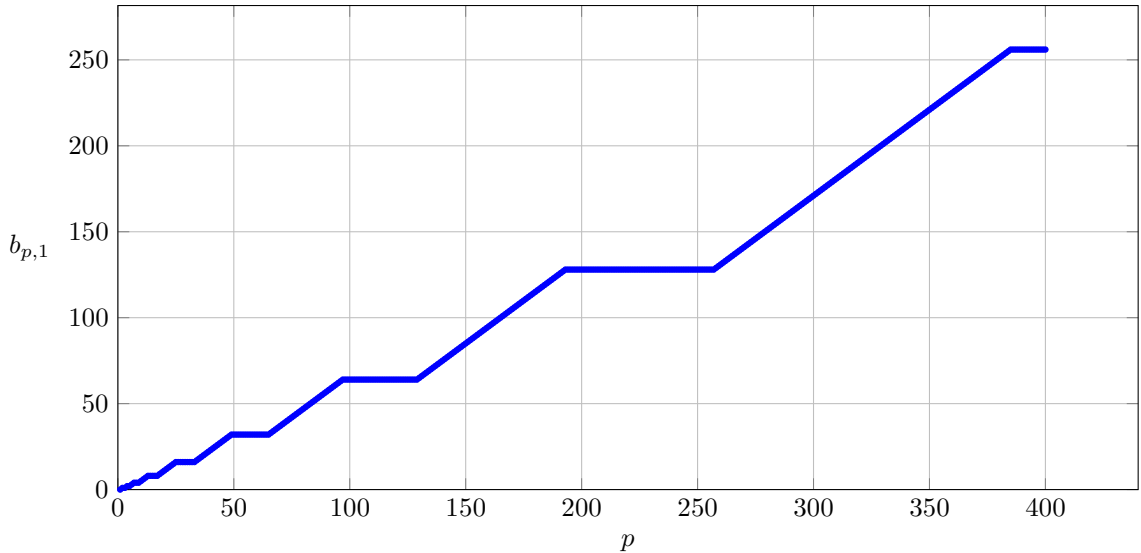
\begin{center}
    \begin{tikzpicture}
        \begin{axis} [xlabel={$p$}
        , ylabel={$b_{p,1}$}
        , ylabel style={rotate=-90}
        , xmin=0
        , ymin = 0
        , grid=major
        , legend pos=north west]
            \addplot[only marks,mark=*, mark size = 1pt, color=blue] table[col sep=semicolon]{b_n,k/b_n,1.csv};
        \end{axis}
    \end{tikzpicture}
\end{center}

\noindent \textbf{Behauptung 3.2:} Mit $m = \left \lfloor \log_{2}(\frac{p-1}{3}) \right \rfloor$ für $p \geq 4$ gilt für alle $p$
\begin{equation*}
    b_{p,1} =
    \begin{cases}
        0 \ \ & \text{für} \ \ p=1 \\
        1 \ \ & \text{für} \ \ p=2 \ \text{oder} \ p=3 \\
        2^{m+1} \ \ & \text{für} \ \ 1+2^{m} \cdot 3 \leq p \leq 1+2^{m} \cdot 4 \\
        p-1-2^{m+1} \ \ & \text{für} \ \ 1+2^{m} \cdot 4 < p < 1+2^{m+1}\cdot 3
    \end{cases}
\end{equation*}

\noindent \textbf{Beweis:} Es wird die Situation betrachtet, dass der Platz links außen bereits durch die erste Person besetzt ist. Zunächst ist $b_{1,1}=0$, da der einzige Platz nach Voraussetzung schon besetzt ist. Ferner ist $b_{2,1}=1$, da der verbliebene leere Platz einen Abstand von $1$ besitzt. Wenn $p=3$, dann besitzt der Platz rechts außen einen Abstand von $2$. Der Platz zwischen den Plätzen links und rechts außen besitzt einen Abstand von $1$, womit $b_{3,1}=1$. \\

\noindent Nun wollen wir eine allgemeine Herleitung der Werte von $b_{p,1}$ im Intervall $[1+2^{m} \cdot 3 ;1+2^{m+1} \cdot 3)$ für alle $m \in \mathbb{N}_0$ durchführen. Da für $p \geq 4$ der Platz rechts außen, der als zweites besetzt wird, mindestens einen Abstand von $3$ besitzt, sei von nun an immer die Situation gegeben, dass neben dem Platz links außen auch der Platz rechts außen fest besetzt ist. Ferner sollen die Schemata aus Abschnitt 2 für das korrekte Hinzufügen von leeren oder besetzten Plätzen verwendet und immer wieder fortgesetzt werden. \\

\noindent Wenn $p=1+2^{m} \cdot 3$, dann gibt es zunächst eine Teilreihe der Länge $2^m \cdot 3 - 1$. Nun kann man diese durch das Besetzen des Platzes in der Mitte dieser Teilreihe in zwei Teilreihen der Länge $\frac{(2^m \cdot 3 - 1)-1}{2} = 2^{m-1} \cdot 3-1$ aufspalten und so immer wieder vorgehen. Man erhält zu bestimmten Zeitpunkten:
\begin{itemize}
    \item zwei Teilreihen der Länge $\frac{(2^m \cdot 3 - 1)-1}{2} = 2^{m-1} \cdot 3-1$
    \item vier Teilreihen der Länge $\frac{(2^{m-1} \cdot 3 - 1)-1}{2} = 2^{m-2} \cdot 3-1$
    \item ...
    \item $2^m$ Teilreihen der Länge $\frac{(2^1 \cdot 3 - 1)-1}{2} = 3-1 = 2$
    
\end{itemize}

\noindent Somit ist dort $b_{p,1}=2^m \cdot 2 = 2^{m+1}$. Nun kann man in die Mitte jeder der $2^m$ Teilreihen der Länge $2$ einen weiteren besetzten Platz nach Schema 2 (mit $S_{2^m}$) hinzufügen, ohne dass sich $b_{p,1}$ verändert. Somit:
\begin{gather}
    b_{p,1} = 2^{m+1} \nonumber \\
    \text{für} \\
    1+2^{m} \cdot 3 \leq p \leq 1+2^{m} \cdot 3 + 2^m = 1+2^{m} \cdot 4 \nonumber
\end{gather}

\noindent Bei $p= 1+2^{m} \cdot 4$ liegen $2^{m+1}$ Teilreihen der Länge $1$ vor. Fügt man nun in eine dieser Teilreihen einen leeren Platz hinzu, dann liegt eine neue Teilreihe der Länge $2$ vor und $b_{p,1}$ steigt um $1$. Nach Schema 1 (mit $S_{2^{m+1}}$) kann man so in jede der $2^{m+1}$ Teilreihen einen leeren Platz hinzufügen. Der Abstand zu $1+2^{m} \cdot 4$ gibt an, zu wie vielen Teilreihen ein leerer Platz hinzugefügt wurde. Es gilt:
\begin{gather}
    b_{p,1} = p-(1+2^{m} \cdot 4)+2^{m+1} = p-1-2^{m+1} \nonumber \\
    \text{für} \\
    1+2^{m} \cdot 4 < p < 1+2^{m} \cdot 4 + 2^{m+1} = 1+2^{m+1} \cdot 3\nonumber
\end{gather}

\noindent Weiter gilt:
\begin{align*}
    p &= 1+2^{m} \cdot 3 \\
    \Leftrightarrow \frac{p-1}{3} &= 2^{m} \\
    \Leftrightarrow m &= \log_{2} \left ( \frac{p-1}{3} \right )
\end{align*}

\noindent Da $m$ im Intervall $[1+2^{m} \cdot 3 ;1+2^{m+1} \cdot 3)$ den gleichen ganzzahligen Wert annehmen soll, muss $m = \left \lfloor \log_{2}(\frac{p-1}{3}) \right \rfloor$. Mit $1+2^0 \cdot 3 = 4$ als kleinste Zahl, die in einem solchen Intervall liegt, sei $m = \lfloor \log_{2}(\frac{p-1}{3}) \rfloor$ für $p \geq 4$. Auf Grund von diesem Aspekt, den Fällen am Anfang, (8) und (9) folgt die Behauptung. \hfill $\square$

\section{Berechnung von $d_{p,k}$}

\noindent Sei $d_{p,k}$ die Anzahl unterschiedlicher Mengen zweier benachbarter Plätze, bei denen im Verlaufe des Besetzens, wenn bei $p \in \mathbb{N}$ Plätzen die erste Person den Platz links außen besetzt, zu einem Zeitpunkt beide Plätze den Abstand $k \in \mathbb{N}$ besitzen.\\

\noindent Damit eine solche Menge existiert, müssen links und rechts von einem besetzten Platz $P$ mit Abstand $k$ eine Teilreihe der Länge $k$ und eine Teilreihe der Länge $k-1$ bzw. vor dem Besetzen von $P$ eine Teilreihe mit gerader Länge vorhanden sein. Setzt man sich nämlich von $P$ auf den benachbarten Platz in der Teilreihe mit der Länge $k$, dann liegen erneut eine Teilreihe der Länge $k$ und eine Teilreihe der Länge $k-1$ vor. In beiden Fällen hat dann der besetzte Platz den Abstand $k$. Die Person auf Platz $P$ konnte also zwischen zwei Möglichkeiten auswählen (siehe dazu auch Abschnitt 1: Teilreihen mit geradem $L$).

\subsection{Berechnung von $d_{p,k}$ für $k \geq 2$}

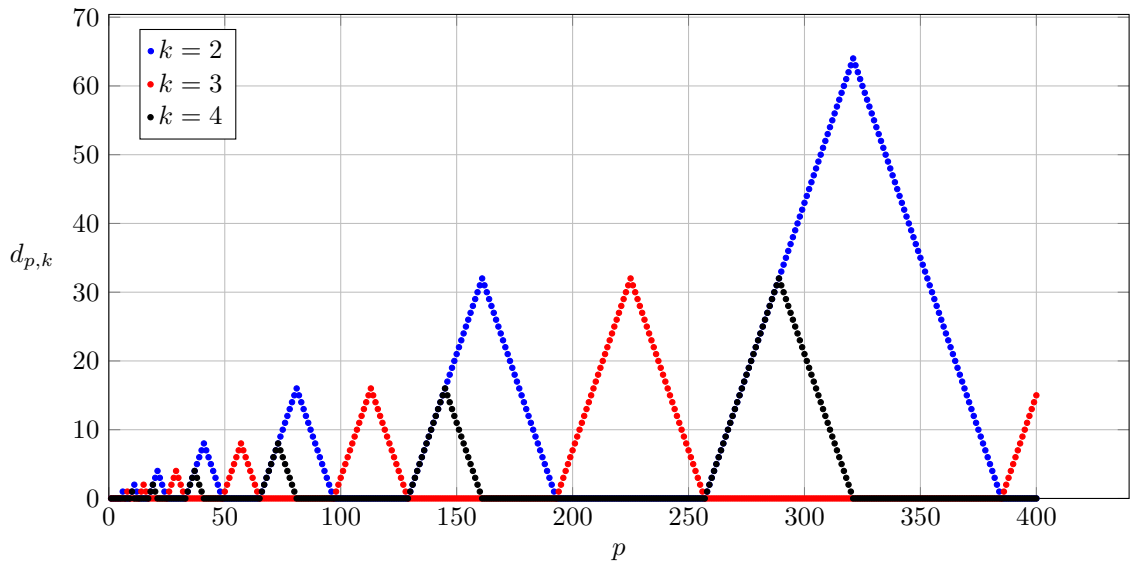
\begin{center}
    \begin{tikzpicture}
        \begin{axis} [xlabel={$p$}
        , ylabel={$d_{p,k}$}
        , ylabel style={rotate=-90}
        , xmin=0
        , ymin = 0
        , grid=major
        , legend pos=north west]
            \addplot[only marks,mark=*, mark size = 1pt, color=blue] table[col sep=semicolon]{d_n,k/d_n,2.csv};
            \addlegendentry{$k=2$}
            \addplot[only marks,mark=*, mark size = 1pt, color=red] table[col sep=semicolon]{d_n,k/d_n,3.csv};;
            \addlegendentry{$k=3$}
            \addplot[only marks,mark=*, mark size = 1pt, color=black] table[col sep=semicolon]{d_n,k/d_n,4.csv};
            \addlegendentry{$k=4$}
        \end{axis}
    \end{tikzpicture}
\end{center}

\noindent \textbf{Behauptung 4.1:} Mit $m = \lfloor \log_{2}(\frac{p-1}{2k}) \rfloor$ für $p \geq 5$ gilt für alle $k \geq 2$ und alle $p$
\begin{equation*}
    d_{p,k} =
    \begin{cases}
        0 
        & \text{für} \ \  p < 1+2k \\

        p - 1-2^{m}\cdot 2k
        & \text{für} \ \ 
        1+2^{m}\cdot 2k \leq p \leq 1+2^{m}(2k+1) \\

        1+2^{m}(2k+2)-p
        & \text{für} \ \
        1+2^{m}(2k+1) < p \leq 1+2^{m}(2k+2) \\

        0 
        & \text{für} \ \
        1+2^{m}(2k+2) < p < 1+2^{m+1} \cdot 2k
    \end{cases}
\end{equation*}

\noindent \textbf{Beweis:} Im Folgenden sei $k\geq 2$. Zudem soll immer die Situation, dass die erste Person den Platz links außen besetzt, betrachtet werden. Da es für $p < 1+2k$ nur bei $p=k+1$ einen Platz mit Abstand $k$ gibt, gilt:
\begin{equation}
    d_{p,k} = 0 \ \ \text{für} \ \  p < 1+2k
\end{equation}

\noindent Analog zu Behauptung 3.1 wollen wir für alle $m \in \mathbb{N}_0$ eine allgemeine Herleitung der Werte von $d_{p,k}$ im Intervall $[1+2^{m}\cdot 2k \ ; \ 1+2^{m+1}\cdot 2k)$ durchführen. Zudem sollen auch hier die Schemata aus Abschnitt 2 für das korrekte Hinzufügen von leeren oder besetzten Plätzen verwendet und immer wieder fortgesetzt werden.\\

\noindent Bei $p=1+2^{m}\cdot 2k$ hat jeder Platz der $2^m$ Plätze mit Abstand $k$ links und rechts von sich eine Teilreihe der Länge $k-1$ (siehe Beweis von Behauptung 3.1), womit $d_{p,k}=0$. Fügt man nun in eine der $2^{m+1}$ Teilreihen einen Platz hinzu, dann hat ein Platz links und rechts von sich eine Teilreihe der Länge $k-1$ und eine Teilreihe der Länge $k$. Da nun eine Menge obiger Definition existiert, erhöht sich $d_{p,k}$ um 1. Nach Schema 1 (mit $S_{2^{m+1}}$) kann man nun insgesamt von $p=1+2^{m}\cdot 2k$ aus $2^m$ leere Plätze hinzufügen und jedes Mal eine neue Menge erzeugen. Somit ist der Abstand zu $1+2^{m}\cdot 2k$ relevant, da dieser angibt, wie viele Mengen durch Hinzufügen von Plätzen erzeugt wurden. Somit:
\begin{gather}
    d_{p,k} = p - (1+2^{m}\cdot 2k) = p - 1-2^{m}\cdot 2k \nonumber \\ 
    \text{für} \\ 
    1+2^{m}\cdot 2k \leq p \leq 1+2^{m}\cdot 2k + 2^m = 1+2^{m}(2k+1) \nonumber
\end{gather}

\noindent Bei $p=1+2^{m}(2k+1)$ hat jeder Platz mit Abstand $k$ links und rechts eine Teilreihe der Länge $k-1$ und eine Teilreihe der Länge $k$. Fügt man dann einen leeren Platz in eine Teilreihe der Länge $k-1$ hinzu, dann wird $d_{p,k}$ um $1$ kleiner, da dann ein neuer Platz links und rechts von sich eine Teilreihe der Länge $k$ besitzt. Damit wird der Abstand zu $p=1+2^m(2k+1)$ relevant, da dieser angibt, zu wie vielen Teilreihen der Länge $k-1$ ein Platz hinzugefügt wurde. Da es $2^m$ Mengen bei $p=1+2^m(2k+1)$  gibt und dort noch $2^m$ Teilreihen der Länge $k-1$ mit Fortsetzung des Schema 1 zu besetzen sind, ist:
\begin{gather}
    d_{p,k} = 2^m -(p - (1+2^{m}(2k+1))) = 1+2^{m}(2k+2)-p  \nonumber \\ 
    \text{für} \\ 
    1+2^{m}(2k+1) < p \leq 1+2^{m}(2k+1) + 2^m = 1+2^{m}(2k+2) \nonumber
\end{gather}

\noindent Nun hat erst bei $p=1+2^{m+1}\cdot 2k$ wieder jeder Platz mit Abstand $k$ links und rechts von sich eine Teilreihe der Länge $k-1$ (siehe Beweis der Behauptung 3.1). Würde man dann einen Platz hinzufügen, hätte man eine Teilreihe der Länge $k$. Dies würde eine Menge erzeugen, da dann ein Platz links und rechts von sich eine Teilreihe der Länge $k$ und eine Teilreihe der Länge $k-1$ besitzt. Im Umkehrschluss bedeutet das, dass für $p < 1+2^{m+1}\cdot 2k$ keine Menge erzeugt wird. Somit ist:
\begin{gather}
    d_{p,k} = 0 \nonumber \\ 
    \text{für} \\ 
    1+2^{m}(2k+2) < p < 1+2^{m+1} \cdot 2k \nonumber
\end{gather}

\noindent Mit (10),(11),...,(13) folgt dann: 
\begin{equation*}
    d_{p,k} =
    \begin{cases}
        0 
        & \text{für} \ \  p < 1+2k \\

        p - 1-2^{m}\cdot 2k
        & \text{für} \ \ 
        1+2^{m}\cdot 2k \leq p \leq 1+2^{m}(2k+1) \\

        1+2^{m}(2k+2)-p
        & \text{für} \ \
        1+2^{m}(2k+1) < p \leq 1+2^{m}(2k+2) \\

        0 
        & \text{für} \ \
        1+2^{m}(2k+2) < p < 1+2^{m+1} \cdot 2k
    \end{cases}
\end{equation*}

\noindent Da $[1;2k] \cup [1+2^{0}\cdot 2k \ ;1+2^{1}\cdot 2k)
\cup [1+2^{1}\cdot 2k \ ;1+2^{2}\cdot 2k)
\cup [1+2^{2}\cdot 2k \ ;1+2^{3}\cdot 2k) 
\cup ... = [1; \infty)$, kann man also mit (10) bis (13) $d_{p,k}$ für alle $p$ berechnen. Weil $m$ im Intervall $[1+2^{m}\cdot 2k \ ;1+2^{m+1}\cdot 2k)$ den gleichen ganzzahligen Wert annehmen soll, ist nach der Herleitung im Beweis von Behauptung 3.1 $m = \lfloor \log_{2}(\frac{p-1}{2k}) \rfloor$. Mit $1+2^0 \cdot 4 = 5$ als kleinste Zahl, die in einem solchen Intervall liegt, sei $m = \lfloor \log_{2}(\frac{p-1}{2k}) \rfloor$ für $p \geq 5$. \hfill $\square$

\subsection{Berechnung von $d_{p,1}$}

\noindent \textbf{Behauptung 4.2:} Mit $m = \lfloor \log_{2}(\frac{p-1}{3}) \rfloor$ für $p \geq 4$ gilt für alle $p$
\begin{equation*}
    d_{p,1} = 
    \begin{cases}
        0 
        & \text{für} \ p < 4 \\
        
        1+2^m \cdot 4-p
        & \text{für} \ 1+2^m \cdot 3 \leq p \leq 1+2^m \cdot 4 \\

        p-1-2^m \cdot 4
        & \text{für} \ 1+2^m \cdot 4 < p < 1+2^{m+1} \cdot 3 \\
    \end{cases}
\end{equation*}

\noindent \textbf{Beweis:} Sei immer die Situation gegeben, dass die erste Person den Platz links außen besetzt. Damit eine Menge entsprechender Definition (siehe Abschnitt 4) existiert, müssen mindestens zwei besetzte Plätze und mindestens zwei leere Plätze vorhanden sein, womit $d_{p,1}$ für $p < 4$. \\

\noindent Nun wollen wir eine allgemeine Herleitung der Werte von $d_{p,1}$ im Intervall $[1+2^{m} \cdot 3 ;1+2^{m+1} \cdot 3)$ für alle $m \in \mathbb{N}_0$ durchführen. Zudem sei von nun an immer die Situation, dass neben dem Platz links außen auch der Platz rechts außen fest besetzt ist, gegeben. Ferner sollen die Schemata aus Abschnitt 2 für das korrekte Hinzufügen von leeren oder besetzten Plätzen verwendet und immer wieder fortgesetzt werden. \\

\noindent Nach Behauptung 3.2 gibt es $2^{m}$ Teilreihen der Länge $2$ bei $p=1+2^m \cdot 3$, womit dort $d_{p,1}=2^{m}$. Fügt man nun nach Schema 2 (mit $S_{2^{m}}$) besetzte Plätze ein, dann sinkt die Anzahl der Mengen pro neu hinzugefügtem besetzten Platz um $1$. Da der Abstand zu $1+2^m \cdot 3$ die Anzahl der neu hinzugefügten Plätze angibt und man insgesamt $2^{m}$ Plätze hinzufügen kann, gilt:
\begin{gather}
    d_{p,1} = 2^{m} - (p-1-2^m \cdot 3) = 1+2^m \cdot 4 - p \nonumber \\
    \text{für} \\
    1+2^m \cdot 3 \leq p \leq 1+2^m \cdot 3+2^m = 1+2^m \cdot 4 \nonumber
\end{gather}

\noindent Bei $p=1+2^m \cdot 4$ liegen also $2^{m+1}$ Teilreihen der Länge $1$ vor. Fügt man in eine dieser Teilreihen einen leeren Platz hinzu, dann entsteht eine neue Menge. Da man so nach Schema 1 (mit $S_{2^{m+1}}$) in alle $2^{m+1}$ Teilreihen der Länge $1$ einen Platz hinzufügen kann, gilt:
\begin{gather}
    d_{p,1} = p -1-2^m \cdot 4 \nonumber \\
    \text{für} \\
    1+2^m \cdot 4 < p < 1+2^m \cdot 4+2^{m+1} = 1+2^{m+1} \cdot 3 \nonumber
\end{gather}

\noindent Da $m$ im Intervall $[1+2^{m} \cdot 3 ;1+2^{m+1} \cdot 3)$ den gleichen ganzzahligen Wert annehmen soll, muss nach dem Beweis von Behauptung 3.2 $m = \left \lfloor \log_{2}(\frac{p-1}{3}) \right \rfloor$. Mit $1+2^0 \cdot 3 = 4$ als kleinste Zahl, die in einem solchen Intervall liegt, sei $m = \lfloor \log_{2}(\frac{p-1}{3}) \rfloor$ für $p \geq 4$. Mit diesem Aspekt, $d_{p,1}=0$ für $p < 4$, (14) und (15) folgt die Behauptung. \hfill $\square$ \\

\section{Berechnung von $a_n$}

\noindent \textbf{Behauptung 5:} Mit $d_{p,1}=0$ für alle $p$ gilt
\begin{equation*}
    a_n = \sum ^{n} _{i=1} \limits \prod \limits ^{n-1} _{j=1} 2^{d_{i,j} + d_{n+1-i,j}} \cdot (b_{i,j} + b_{n+1-i,j})!
\end{equation*}

\noindent \textbf{Beweis:} Zunächst kann die erste Person den $i$-ten Platz ($i \in \mathbb{N} \ \text{und} \ i \in [1;n]$) von links besetzen. Diese Fälle müssen einzeln betrachtet und die jeweiligen Möglichkeiten hinterher addiert werden. Zudem kann jeder solcher Fall als Paar $(n;i)$ dargestellt werden. Besetzt nun die erste Person den $i$-ten Platz von links, dann gibt es $i-1$ leere Plätze links und $n-i$ leere Plätze rechts von diesem Platz. Nun kann man die $i-1$ leeren Plätze links vom Platz der ersten Person und den Platz der ersten Person als neue Situation $(i;1)$ auffassen. In ähnlicher Weise kommt man mit den $n-i$ leeren Plätzen rechts vom besetzten Platz zur neuen Situation $(n+1-i;1)$. \\

\noindent Setzt sich die erste Person nun auf den $i$-ten Platz von links, dann gibt es also $b_{i,1}+b_{n+1-i,1}$ Plätze mit Abstand $1$, $b_{i,2}+b_{n+1-i,2}$ Plätze mit Abstand $2$, ..., $b_{i,n-1}+b_{n+1-i,n-1}$ Plätze mit dem größtmöglichen Abstand $n-1$. Nun gibt es also bei Abstand $j \in \mathbb{N}$ $(b_{i,j}+b_{n+1-i,j})!$ Möglichkeiten, diese $b_{i,j}+b_{n+1-i,j}$ Plätze mit $b_{i,j}+b_{n+1-i,j}$ Personen zu besetzen. Zudem muss vor dem Besetzen von Plätzen mit Abstand $j$ festgelegt werden, welcher Platz einer Menge (siehe Abschnitt 4) besetzt werden darf. Bei Abstand $j$ gibt es dafür $2^{d_{i,j} + d_{n+1-i,j}}$ Möglichkeiten. Ferner sei $d_{p,1}=0$ für alle $p$, da bei einer Teilreihe der Länge $2$ der Aspekt der zwei Möglickeiten wegfällt (siehe Abschnitt 1). Gemäß der Argumentation in Abschnitt 1 werden erst die Plätze mit Abstand $n-1$ und dann die Plätze mit Abstand $n-2/n-3/.../2/1$ besetzt. Folglich erhält man mit der Produktregel der Kombinatorik
\begin{equation*}
    \prod \limits ^{n-1} _{j=1} 2^{d_{i,j} + d_{n+1-i,j}} \cdot (b_{i,j} + b_{n+1-i,j})!
\end{equation*}
\noindent Möglichkeiten, wenn sich die erste Person auf den $i$-ten Platz von links setzt. Mit $i\in [1;n]$ erhält man schlussendlich:
\begin{equation*}
    a_n = \sum ^{n} _{i=1} \limits \prod \limits ^{n-1} _{j=1} 2^{d_{i,j} + d_{n+1-i,j}} \cdot (b_{i,j} + b_{n+1-i,j})!
\end{equation*}

\hfill $\square$

\section{Abschätzungen}
Zunächst ist $n \leq a_n$, da, wenn sich die erste Person auf den ersten/zweiten/.../$n$-ten Platz von links setzt, immer mindestens eine neue Möglichkeit entsteht. Weiter ist $a_n \leq n!$ , da ohne die Höflichkeitsregel es für die erste Person $n$ mögliche Plätze, für die zweite Person $n-1$ mögliche Plätze,..., für die letzte Person einen möglichen Platz gibt.

\subsection{Abschätzung von $b_{p,k}$ für $k \geq 2$}

\noindent \textbf{Lemma 6.1:} Es gilt mit $m' = \lfloor \log_{2}(\frac{p-1}{4}) \rfloor$ 
\begin{equation*}
    \prod _{j=1} ^{m'} \limits (2^{j-1})! 
    \leq
    \prod _{j=2} ^{p-1} \limits (b_{p,j})!
\end{equation*}

\noindent \textbf{Beweis:} Mit $m = \lfloor \log_{2}(\frac{p-1}{2k}) \rfloor$ für $p \geq 5$ gilt nach Behauptung 3.1 für $k \geq 2$
\begin{gather}
    b_{p,k}= 2^m \nonumber \\
    \text{für} \\
    1+2^{m}\cdot 2k \leq p \leq 1+2^{m}(2k+1) \nonumber
\end{gather}
\begin{gather}
    b_{p,k}=1+2^m(2k+2)-p \nonumber \\ 
    \text{für} \\
    1+2^{m}(2k+1) < p \leq 1+2^{m}(2k+2) = 1+2^{m}\cdot 2(k+1) \nonumber
\end{gather}
\begin{gather}
    b_{p,k+1}=0 \nonumber \\ 
    \text{für} \\
    1+2^{m-1}(2(k+1)+2) = 1+2^{m}(k+2) < p \leq 1+2^{m-1}(4(k+1)-2) = 1+2^{m}(2k+1) \nonumber
\end{gather}
\begin{gather}
     b_{p,k+1}=p-1-2^{m-1}(4(k+1)-2) \nonumber \\ 
     \text{für} \\
     1+2^{m-1}(4(k+1)-2) = 1+2^{m}(2k+1) < p < 1+2^{m-1+1}2(k+1) =  1+2^{m}\cdot 2(k+1)  \nonumber
\end{gather}

\noindent Mit (16) und (18) gilt im Intervall $[1+2^{m}\cdot 2k; 1+2^{m}(2k+1)]$ (Bei $k=2$ gibt es den in (16) nicht enthaltenen Fall $p=1+2^m \cdot 2k$. Nach Behauptung 3.1 gilt dort allerdings auch $b_{p,k+1}=0$)
\begin{align*}
    b_{p,k} + b_{p,k+1} = 2^m+0 = 2^m .
\end{align*}
\noindent Mit (17) und (19) gilt im Intervall $(1+2^{m}(2k+1);1+2^{m}\cdot 2 (k+1))$
\begin{align*}
    &b_{p,k} + b_{p,k+1} \\
    &= 1+2^m(2k+2)-p+p-1-2^{m-1}(4(k+1)-2) \\
    &= 2^m(2k+2)-2^{m}(2k+2-1) \\
    &= 2^m(2k+1-2k)=2^{m} .
\end{align*}

\noindent Somit gilt im Intervall $[1+2^m \cdot 2k \ ; \ 1+2^m \cdot2(k+1))$ 
\begin{equation*}
    b_{p,k}+b_{p,k+1}=2^m ,
\end{equation*}

\noindent womit dort immer entweder $b_{p,k}\geq 2^{m-1}$ oder $b_{p,k+1}\geq 2^{m-1}$. Zudem gibt es auf Grund von 
$[1+2^m \cdot 4 \ ; \ 1+2^m \cdot 6) 
\cup [1+2^m \cdot 6 \ ; \ 1+2^m \cdot 8) 
\cup ... = [1+2^m \cdot 4 ; \infty) $ für $p \geq 1+2^m \cdot 4$ (hierbei soll $m$ eine beliebige nichtnegative ganze Zahl sein) immer ein $k$, sodass $b_{p,k}+b_{p,k+1}=2^m$. Nun kann man folgendermaßen berechnen, welche Werte $m$ annimmt:
\begin{align*}
    p &\geq 1+2^m \cdot 4 \\
    \Leftrightarrow \frac{p-1}{4} &\geq 2^m \\
    \Leftrightarrow \log_{2} \left (\frac{p-1}{4} \right ) &\geq m
\end{align*}

\noindent Damit ist $m' = \left \lfloor\log_2 (\frac{p-1}{4}) \right \rfloor$ der größtmögliche Wert von $m$. Mit entweder $b_{p,k}\geq 2^{m-1}$  oder $b_{p,k+1}\geq 2^{m-1}$ für bestimmte $k$ und $p-1$ als größtmöglichem Abstand (Es ist $b_{p,k}=0$ für $k \geq p$) folgt dann:

\begin{equation*}
     \prod _{j=1} ^{m'} \limits (2^{j-1})! 
     \leq
     \prod _{j=2} ^{p-1} \limits (b_{p,j})! 
\end{equation*}

\hfill $\square$

\noindent \textbf{Behauptung 6.1:} Es gilt mit $m_i = \left \lfloor\log_2 (\frac{i-1}{4}) \right \rfloor$ für $n \geq 2$

\begin{equation*}
    2 \cdot \sum 
    ^{n} 
    _{i=\left \lceil \frac{n}{2} \right \rceil +1} \limits
    \prod \limits ^{m_i} _{j=1} (2^{j-1})!  
    \leq
    \sum ^{n} _{i=1} \limits \prod \limits ^{n-1} _{j=2} (b_{i,j} + b_{n+1-i,j})!   
\end{equation*}

\noindent \textbf{Beweis:} Es ist für $n \geq 2$:
\begin{itemize}
    \item $b_{n,k}=b_{n+1-1,k} \ (i=n \ \text{und} \ i=1 \ \text{bei Formel für} \ a_n)$
    \item $b_{n-1,k}=b_{n+1-2,k} \ (i=n-1 \ \text{und} \ i=2 \ \text{bei Formel für} \ a_n)$
    \item $b_{n-2,k}=b_{n+1-3,k} \ (i=n-2 \ \text{und} \ i=3 \ \text{bei Formel für} \ a_n)$
    \item ...
    \item $b_{\left \lceil \frac{n}{2} \right \rceil +1,k}
    = b_{n+1-\left \lfloor \frac{n}{2} \right \rfloor,k}
    \ 
    (i=\left \lceil \frac{n}{2} \right \rceil +1 \ \text{und} 
    \ i=\left \lfloor \frac{n}{2} \right \rfloor
    \ \text{bei Formel für} \ a_n)$
\end{itemize}

\noindent Nun kann man damit nach unten abschätzen:

\begin{equation*}
    2 \cdot \sum 
    ^{n} 
    _{i=\left \lceil \frac{n}{2} \right \rceil +1} \limits 
    \prod \limits ^{i-1} _{j=2} (b_{i,j})!
    \leq
    \sum ^{n} _{i=1} \limits \prod \limits ^{n-1} _{j=2} (b_{i,j} + b_{n+1-i,j})!
\end{equation*}

\noindent Mit Lemma 6.1 ist dann mit $m_i = \left \lfloor\log_2 (\frac{i-1}{4}) \right \rfloor$:

\begin{equation*}
    2 \cdot \sum 
    ^{n} 
    _{i=\left \lceil \frac{n}{2} \right \rceil +1} \limits 
    \prod \limits ^{m_i} _{j=1} (2^{j-1})!
    \leq
    2 \cdot \sum 
    ^{n} 
    _{i=\left \lceil \frac{n}{2} \right \rceil +1} \limits 
    \prod \limits ^{i-1} _{j=2} (b_{i,j})!
\end{equation*}

\hfill $\square$

\subsection{Abschätzung von $b_{p,1}$}

\noindent \textbf{Behauptung 6.2.1:} Es gilt für alle $p$
\begin{equation*}
    \left \lfloor \frac{p-1}{2} \right \rfloor
    \leq
    b_{p,1}
\end{equation*}

\noindent \textbf{Beweis:} $b_{p,1}$ hat nach Behauptung 3.2 im Intervall $[1+2^{m}\cdot 3 \ ;1+2^{m+1}\cdot 3)$ den minimalen Wert $2^{m+1}$. Dieser Wert tritt im genannten Intervall zuletzt bei $p=1+2^{m}\cdot 4$ auf. Nun definiert man sich zwei Punkte $P$ und $Q$ ($c \in \mathbb{N}$): 
\begin{equation*}
    P(1+2^{m}\cdot 4 \mid 2^{m+1}) \ \ \ \ Q(1+2^{m+c}\cdot 4 \mid 2^{m+1+c})
\end{equation*}

\noindent Betrachtet man nun den Differenzenquotient, dann stellt man fest, dass alle Punkte dieser Art auf einer Geraden liegen: \\
\begin{equation*}
    m_{\text{Steigung}} = \frac{2^{m+1+c}-2^{m+1}}{1+2^{m+c}\cdot 4-(1+2^{m}\cdot 4)} = \frac{2^{m+1} \cdot (2^c - 1)}{2^{m} \cdot 4 \cdot (2^c-1)} = \frac{2}{4} = \frac{1}{2}
\end{equation*}

\noindent Einsetzen der Koordinaten von $P$ und $m_{\text{Steigung}}$ in die Funktionsgleichung einer Gerade $g$ mit $x \in \mathbb{R}$ ergibt:
\begin{align*}
    g(x) &= m_{\text{Steigung}} \cdot x + f \\
    2^{m+1} &= \frac{1}{2} \cdot (1+2^{m}\cdot 4)+f \\
    2^{m+1} &= \frac{1}{2} + 2^{m+1}+f \\
    \Leftrightarrow f &= -\frac{1}{2}
\end{align*}
\begin{equation*}
    \Rightarrow g(x)= \frac{1}{2} \cdot x -\frac{1}{2} = \frac{1}{2} \cdot (x-1)
\end{equation*}

\noindent Abschließend gilt für $p \geq 1+2^0 \cdot 3 = 4$
\begin{align*}
    \left \lfloor \frac{p-1}{2} \right \rfloor
    = \left \lfloor \frac{1}{2} \cdot (p-1) \right \rfloor 
    = \left \lfloor g(p) \right \rfloor
    \leq g(p)
    \leq b_{p,1} \ .
\end{align*}

\noindent Da zudem 
$\left \lfloor \frac{1-1}{2} \right \rfloor = 0 \leq b_{1,1} = 0 ,$
$\left \lfloor \frac{2-1}{2} \right \rfloor = 0 \leq b_{2,1} = 1 $ und
$\left \lfloor \frac{3-1}{2} \right \rfloor = 1 \leq b_{3,1} = 1 $ gilt 
$\left \lfloor \frac{p-1}{2} \right \rfloor \leq b_{p,1} $ für alle $p$.

\hfill $\square$

\noindent \textbf{Behauptung 6.2.2:} Es gilt für alle $p$
\begin{equation*}
    b_{p,1}
    \leq
    \left \lceil \frac{2\cdot (p-1)}{3} \right \rceil
\end{equation*}

\noindent \textbf{Beweis:} $b_{p,1}$ hat nach Behauptung 3.2 im Intervall $(1+2^{m}\cdot 3 \ ;1+2^{m+1}\cdot 3]$ den maximalen Wert $2^{m+2}$. Dieser Wert tritt im genannten Intervall zuerst bei $n=1+2^{m+1}\cdot 3$ auf. Nun definieren wir uns zwei allgemeine Punkte $P$ und $Q$ ($c \in \mathbb{N}$):
\begin{align*}
    P(1+2^{m+1}\cdot 3 \mid 2^{m+2}) \ \ \ \ \ Q(1+2^{m+1+c}\cdot 3 \mid 2^{m+2+c})
\end{align*}

\noindent Betrachtet man nun den Differenzenquotient, dann stellt man fest, dass alle Punkte dieser Art auf einer Geraden liegen: \\
\begin{equation*}
    m_{\text{Steigung}} = \frac{2^{m+2+c}-2^{m+2}}{1+2^{m+1+c}\cdot 3-(1+2^{m+1}\cdot 3)} = \frac{2^{m+2} \cdot (2^c - 1)}{2^{m+1} \cdot (2^c-1) \cdot 3} = \frac{2}{3}
\end{equation*}

\noindent Einsetzen der Koordinaten von $P$ und $m_{\text{Steigung}}$ in die Funktionsgleichung einer Gerade $g$ mit $x \in \mathbb{R}$:
\begin{align*}
    g(x) &= m_{\text{Steigung}} \cdot x + f \\
    2^{m+2} &= \frac{2}{3} \cdot (1+2^{m+1}\cdot 3)+f \\
    2^{m+2} &= \frac{2}{3} + 2^{m+2} + f \\
    \Leftrightarrow f &= -\frac{2}{3}
\end{align*}
\begin{equation*}
    \Rightarrow g(x)= \frac{2}{3} \cdot x -\frac{2}{3} = \frac{2(x-1)}{3}
\end{equation*}

\noindent Schließlich ist für $p \geq 1+2^0 \cdot 3 + 1 =5$:
\begin{align*}
    b_{p,1} 
    \leq g(p) 
    \leq \lceil g(p) \rceil 
    = \left \lceil \frac{2(p-1)}{3} \right \rceil
\end{align*}

\noindent Des Weiteren gilt $b_{1,1} = 0 \leq \left \lceil \frac{2(1-1)}{3} \right \rceil = 0, 
b_{2,1} = 1 \leq \left \lceil \frac{2(2-1)}{3} \right \rceil = 1,
b_{3,1} = 1 \leq \left \lceil \frac{2(3-1)}{3} \right \rceil = 2 \ \text{und} \
b_{4,1} = 2 \leq \left \lceil \frac{2(4-1)}{3} \right \rceil = 2$
, womit $b_{p,1} \leq \left \lceil \frac{2(p-1)}{3} \right \rceil$ für alle $p$. \hfill $\square$

\subsection{Abschätzung von $a_n$}

\noindent \textbf{Behauptung 6.3.1:} Es gilt mit $m_i = \left \lfloor\log_2 (\frac{i-1}{4}) \right \rfloor$ für $n \geq 2$
\begin{align*}
    \underbrace{
    2 \cdot \sum 
    ^{n} 
    _{i=\left \lceil \frac{n}{2} \right \rceil +1} \limits 
    \left (\left \lfloor \frac{i-1}{2} \right \rfloor
    + \left \lfloor \frac{n-i}{2} \right \rfloor \right )! 
    \cdot
    \prod \limits ^{m_i} _{j=1} (2^{j-1})!
    } _{= \ U}
    \leq
    a_n
\end{align*}

\noindent \textbf{Beweis:} Zunächst gilt mit $ 0 \leq d_{i,j}$ bzw. $ 0 \leq d_{n+1-i,j}$
\begin{align*}
    \sum ^{n} _{i=1} \limits 
    (b_{i,1}+b_{n+1-i,1})! \cdot
    \prod \limits ^{n-1} _{j=2} (b_{i,j} + b_{n+1-i,j})!
    \leq 
    \sum ^{n} _{i=1} \limits \prod \limits ^{n-1} _{j=1} 2^{d_{i,j} + d_{n+1-i,j}} \cdot (b_{i,j} + b_{n+1-i,j})!
    = a_n
\end{align*}

\noindent Mit Behauptung 6.1 gilt dann mit $m_i = \left \lfloor\log_2 (\frac{i-1}{4}) \right \rfloor$ für $n \geq 2$:
\begin{align*}
    2 \cdot \sum 
    ^{n} 
    _{i=\left \lceil \frac{n}{2} \right \rceil +1} \limits
    (b_{i,1}+b_{n+1-i,1})! \cdot
    \prod \limits ^{m_i} _{j=1} (2^{j-1})! 
    \leq
    \sum ^{n} _{i=1} \limits 
    (b_{i,1}+b_{n+1-i,1})! \cdot
    \prod \limits ^{n-1} _{j=2} (b_{i,j} + b_{n+1-i,j})!
\end{align*}

\noindent Mit Behauptung 6.2.1 folgt schließlich:
\begin{align*}
    2 \cdot \sum 
    ^{n} 
    _{i=\left \lceil \frac{n}{2} \right \rceil +1} \limits 
    \left (\left \lfloor \frac{i-1}{2} \right \rfloor
    + \left \lfloor \frac{n-i}{2} \right \rfloor \right )! 
    \cdot
    \prod \limits ^{m_i} _{j=1} (2^{j-1})!
    \leq
    2 \cdot \sum 
    ^{n} 
    _{i=\left \lceil \frac{n}{2} \right \rceil +1} \limits 
    (b_{i,1}+b_{n+1-i,1})! \cdot
    \prod \limits ^{m_i} _{j=1} (2^{j-1})! 
\end{align*} \hfill $\square$

\noindent \textbf{Behauptung 6.3.2:} Es gilt für alle $n$
\begin{align*}
    a_n 
    \leq 
    \underbrace{
    \sum _{i=1} ^n \limits
    \left ( \left \lceil \frac{2(i-1)}{3} \right \rceil +\left \lceil \frac{2(n-i)}{3} \right \rceil \right )! \cdot
    \left ( n-\left \lceil \frac{2(i-1)}{3} \right \rceil-\left \lceil \frac{2(n-i)}{3} \right \rceil \right )!
    } _{= \ O}
\end{align*}

\noindent \textbf{Beweis:} Setzt sich die erste Person auf den $i$-ten Platz von links ($i \in [1;n]$), dann besetzen $b_{i,1}+b_{n+1-i,1}$ Personen einen Platz mit Abstand $1$. Nun gibt es $(b_{i,1}+b_{n+1-i,1})!$ Möglichkeiten, wie die $b_{i,1}+b_{n+1-i,1}$ Personen diese $b_{i,1}+b_{n+1-i,1}$ Plätze besetzen können. Ferner gibt es damit für die erste Person höchstens $n-b_{i,1}-b_{n+1-i,1}$ mögliche Plätze, für die zweite Person höchstens $n-b_{i,1}-b_{n+1-i,1}-1$ mögliche Plätze, ..., für die $n-b_{i,1}-b_{n+1-i,1}$-te Person höchstens einen möglichen Platz. Somit folgt:
\begin{align*}
    a_n \leq \sum _{i=1} ^n \limits
    (b_{i,1}+b_{n+1-i,1})! \cdot
    (n-b_{i,1}-b_{n+1-i,1})!
\end{align*}

\noindent Mit Behauptung 6.2.2 ist weiter:
\begin{align*}
    & \sum _{i=1} ^n \limits
    (b_{i,1}+b_{n+1-i,1})! \cdot
    (n-b_{i,1}-b_{n+1-i,1})! \\
    & \leq 
    \sum _{i=1} ^n \limits
    \left ( \left \lceil \frac{2(i-1)}{3} \right \rceil +\left \lceil \frac{2(n-i)}{3} \right \rceil \right )! \cdot
    \left ( n-\left \lceil \frac{2(i-1)}{3} \right \rceil-\left \lceil \frac{2(n-i)}{3} \right \rceil \right )!
\end{align*} \hfill $\square$

\section{Vergleich von Formel mit Abschätzungen}

\begin{center}
    \begin{tabular}{|c|c|c|c|c|c|c|}
    \hline $n$ & $\displaystyle \frac{U}{a_n}$ & $U$ & $a_n$ &  $O$ & $\displaystyle \frac{O}{a_n}$ & $n!$ \\ \hline 
    1 & / & / &  1  & 1 & 1 & 1 \\ \hline
    2 & 1 & 2 &  2  & 2 & 1 & 2 \\ \hline
    3 & 0,5 & 2 &  4  & 6 & 1,5 & 6\\ \hline
    4 & 0,5 & 4 &  8  & 20 & 2,5 & 24\\ \hline
    5 & 0,3 & 6 &  20  & 72 & 3,6 & 120\\ \hline
    6 & 0,25 & 12 &  48  & 288 & 6 & 720\\ \hline
    7 & $\approx 0,1296$  & 28 & 216  & 1392 & $\approx 6,44$ & 5040\\ \hline
    8 & $\approx 0,0833$ & 48 &  576  & 7200 & 12,5 & 40320\\ \hline
    9 & $\approx 0,0862$ & 120 &  1392  & 38880 & $\approx 27,93$ & 362880\\ \hline
    10 & $\approx 0,0333$ &240 &  7200  & 250560 & 34,8 & 3628800\\ \hline
    15 & $\approx 0,0021$ & 44640 &  21611520 & 6531840000 & $\approx 302,24$  & 1307674368000\\ \hline
    \end{tabular}
\end{center}

\section{Formeln für bereits vorhandene OEIS-Folgen}

\noindent Sei $A166079(n)$ das $n$-te Folgenglied von $A166079$ (siehe [2]). Ähnlich soll mit anderen Folgen wie z.B. A095236 (siehe [3]) verfahren werden.

\subsection{A166079}

\noindent Zunächst kann man mit den Erkenntnissen aus den Abschnitten 3 und 5 eine Formel für die Situation in [2] angeben:
\begin{equation*}
    A166079(n) = n - (b_{1,1}+b_{n,1}) = n - b_{n,1}
\end{equation*}

\subsection{A095236}

\noindent Die in [3] dargestellte Situation unterscheidet sich derart von der in Abschnitt 1 beschriebenen Situation, dass eine unterschiedliche Höflichkeitsregel zum Besetzen der Plätze vorliegt: \\

\noindent \textbf{Höflichkeitsregel:} Wer dazu kommt, wählt einen der Plätze so aus, dass der Abstand zu einem besetzten Platz maximal wird. Dabei soll die jeweilige Person immer einen der Plätze in den Teilreihen mit der größten Länge besetzen. \\

\noindent \textbf{Behauptung 8.2:} Es gilt für alle $n$
\begin{equation*}
    A095236(n) = 
    \sum _{i=1} ^{n} \limits
    \prod _{j=1} ^{n-1} \limits
    2^{ d_{i,j}+d_{n+1-i,j} } \cdot
    (d_{i,j}+d_{n+1-i,j})! \cdot 
    (b_{i,j}+b_{n+1-i,j}-d_{i,j}-d_{n+1-i,j})!
\end{equation*}

\noindent \textbf{Beweis:} Analog zum Beweis von Behauptung 5 kann die Situation, dass bei $n$ Plätzen die erste Person den $i$-ten Platz von links ($i \in \mathbb{N}$ und $i \in [1;n]$) besetzt, als Paar $(n;i)$ dargestellt werden. Dabei kann $(n;i)$ erneut in $(i;1)$ und $(n+1-i;1)$ aufgeteilt werden.\\

\noindent Betrachtet man nun die neue Höflichkeitsregel, dann stellt man beim Besetzen von Plätzen mit Abstand $j \in \mathbb{N}$ fest, dass erst von jeder Menge ein Platz besetzt wird, bevor die anderen Plätzen mit Abstand $j$ besetzt werden. Da beim Besetzen der Plätzen der Mengen beide Plätze einer Menge in Frage kommen und es $d_{i,j}+d_{n+1-i,j}$ Mengen gibt, gibt es also $2^{d_{i,j}+d_{n+1-i,j}} \cdot (d_{i,j}+d_{n+1-i,j})!$ Möglichkeiten, erst von jeder Menge einen Platz zu besetzen. Für $b_{i,j}+b_{n+1-i,j}-d_{i,j}-d_{n+1-i,j}$ als Anzahl der verliebenen Plätze mit Abstand $j$ gibt es also $(b_{i,j}+b_{n+1-i,j-}-d_{i,j}-d_{n+1-i,j})!$ Möglichkeiten, diese zu besetzen. Da erst die Plätze mit dem größtmöglichen Abstand $n-1$ und dann die Plätze mit Abstand $n-2/n-3/.../2/1$ besetzt werden, erhält man
\begin{equation*}
    \prod _{j=1} ^{n-1} \limits
    2^{ d_{i,j}+d_{n+1-i,j} } \cdot
    (d_{i,j}+d_{n+1-i,j})! \cdot 
    (b_{i,j}+b_{n+1-i,j}-d_{i,j}-d_{n+1-i,j})!
\end{equation*}
Möglichkeiten, wenn sich die erste Person auf den $i$-ten Platz von links setzt. Mit $i \in [1;n]$ folgt:
\begin{equation*}
    A095236(n) = 
    \sum _{i=1} ^{n} \limits
    \prod _{j=1} ^{n-1} \limits
    2^{ d_{i,j}+d_{n+1-i,j} } \cdot
    (d_{i,j}+d_{n+1-i,j})! \cdot 
    (b_{i,j}+b_{n+1-i,j}-d_{i,j}-d_{n+1-i,j})! \ \ \ \ \ \ \ \ \ \ \ \ \ \ \ \ \ \ \ \ \ \ \hfill \square
\end{equation*} 

\subsection{A095240}

Da [4] die Fälle $i=1$ und $i=n$ der Formel für $A095236(n)$ (siehe Behauptung 8.2) behandelt, kann für $A095240(n)$ auch eine Formel für $n \geq 2$ angegeben werden:
\begin{align*}
    A095240(n) = 2 \cdot
    \prod _{j=1} ^{n-1} \limits
    2^{ d_{1,j}+d_{n+1-1,j} } \cdot
    (d_{1,j}+d_{n+1-1,j})! \cdot 
    (b_{1,j}+b_{n+1-1,j}-d_{1,j}-d_{n+1-1,j})!
\end{align*}

\noindent Mit $d_{1,j}=0$ und $b_{1,j}=0$ für alle $j$ (siehe Behauptungen 3.1, 3.2, 4.1 und 4.2) gilt:
\begin{align*}
    A095240(n) 
    = 2 \cdot
    \prod _{j=1} ^{n-1} \limits
    2^{ d_{n,j} } \cdot
    (d_{n,j})! \cdot 
    (b_{n,j}-d_{n,j})! 
\end{align*}

\subsection{A095912}

\noindent Die Situation in [5] ist eine Variante von der Situation in Abschnitt 8.2. Damit die in [5] beschriebene zusätzliche Regel eine Rolle spielt, muss mindestens eine Reihe der Form $OXOX$ vorhanden sein. Damit kann für $n < 4$ die Formel für $A095236(n)$ aus Behauptung 8.2 verwendet werden. \\

\noindent \textbf{Behauptung:} Es gilt für $n \geq 4$:
\begin{align*}
    A095912(n)
     = & \ 2 \cdot
     \sum _{i=n-1} ^{n} \limits
    \prod _{j=1} ^{n-1} \limits
    2^{ d_{i,j} } \cdot
    (d_{i,j})! \cdot 
    (b_{i,j}-d_{i,j})! \\
    & + \sum _{i=3} ^{n-2} \limits
    \prod _{j=1} ^{n-1} \limits
    2^{ d_{i,j}+d_{n+1-i,j} } \cdot
    (d_{i,j}+d_{n+1-i,j})! \cdot 
    (b_{i,j}+b_{n+1-i,j}-d_{i,j}-d_{n+1-i,j})!
\end{align*}

\noindent \textbf{Beweis:} Die in [5] beschriebene zusätzliche Regel bezieht sich nur auf das Besetzen von Plätzen mit Abstand $1$. Da nur bei $i=2$ und $i=n-1$ ein entsprechender Platz links oder rechts außen vorhanden ist, sei zunächst $i=2$ oder $i=n-1$.\\

\noindent Beim Besetzen von Plätzen mit Abstand $1$ werden nun erst von jeder Menge ein Platz besetzt, wofür es $2^{d_{2,1}+d_{n+1-2,1}} \cdot (d_{2,1}+d_{n+1-2,1})!$ Möglichkeiten gibt (siehe Beweis von Behauptung 8.2). Danach wird der Platz links oder rechts außen besetzt. Anschließend werden die  verbliebenen Plätze mit Abstand $1$ ($b_{2,1}+b_{n+1-2,1}-d_{2,1}-d_{n+1-2,1}-1$ als Anzahl der verblieben Plätze mit Abstand $1$) besetzt. Dafür gibt es $ (b_{2,1}+b_{n+1-2,1}-d_{2,1}-d_{n+1-2,1}-1)!$ Möglichkeiten.\\

\noindent Bei $i=2$ oder $i=n-1$ gibt es somit $2^{d_{2,1}+d_{n+1-2,1}} \cdot (d_{2,1}+d_{n+1-2,1})! \cdot 1 \cdot (b_{2,1}+b_{n+1-2,1}-d_{2,1}-d_{n+1-2,1}-1)!$ Möglichkeiten, die Plätze mit Abstand $1$ zu besetzen. Mit den gleichbleibenden Möglichkeiten bei den anderen Werten von $i$ und den gleichbleibenden Möglichkeiten bei Abständen $\geq 2$ (siehe dazu Abschnitt 8.2) folgt dann für $n \geq 4$:

\begin{align*}
    A095912(n) = \ &2 \cdot
    \prod _{j=1} ^{n-1} \limits
    2^{ d_{1,j}+d_{n+1-1,j} } \cdot
    (d_{1,j}+d_{n+1-1,j})! \cdot 
    (b_{1,j}+b_{n+1-1,j}-d_{1,j}-d_{n+1-1,j})! \\
    & + 2 \cdot 2^{d_{2,1}+d_{n+1-2,1}} \cdot (d_{2,1}+d_{n+1-2,1})! \cdot 1 \cdot (b_{2,1}+b_{n+1-2,1}-d_{2,1}-d_{n+1-2,1}-1)! \ \cdot \\
    & \prod _{j=2} ^{n-1} \limits
    2^{ d_{2,j}+d_{n+1-2,j} } \cdot
    (d_{2,j}+d_{n+1-2,j})! \cdot
    (b_{2,j}+b_{n+1-2,j}-d_{2,j}-d_{n+1-2,j})! \\
    & + \sum _{i=3} ^{n-2} \limits
    \prod _{j=1} ^{n-1} \limits
    2^{ d_{i,j}+d_{n+1-i,j} } \cdot
    (d_{i,j}+d_{n+1-i,j})! \cdot 
    (b_{i,j}+b_{n+1-i,j}-d_{i,j}-d_{n+1-i,j})!
\end{align*}

\noindent Mit $d_{1,j}=0, d_{2,j}=0 \ \text{und} \ b_{1,j}=0$ für alle $j$, $b_{2,j}=0$ für $j \geq 2$ und $b_{2,1}=1$ (siehe Behauptungen 3.1, 3.2, 4.1 und 4.2) gilt weiter:

\begin{align*}
    A095912(n) = \ &2 \cdot
    \prod _{j=1} ^{n-1} \limits
    2^{ d_{n,j} } \cdot
    (d_{n,j})! \cdot 
    (b_{n,j}-d_{n,j})! \\
    & + 2 \cdot 2^{d_{n-1,1}} \cdot (d_{n-1,1})! \cdot (b_{n-1,1}-d_{n-1,1})! \ \cdot
    \prod _{j=2} ^{n-1} \limits
    2^{d_{n-1,j} } \cdot
    (d_{n-1,j})! \cdot
    (b_{n-1,j}-d_{n-1,j})! \\
    & + \sum _{i=3} ^{n-2} \limits
    \prod _{j=1} ^{n-1} \limits
    2^{ d_{i,j}+d_{n+1-i,j} } \cdot
    (d_{i,j}+d_{n+1-i,j})! \cdot 
    (b_{i,j}+b_{n+1-i,j}-d_{i,j}-d_{n+1-i,j})!
\end{align*}

\noindent Damit gilt auch:
\begin{align*}
    A095912(n) = \ &2 \cdot
    \prod _{j=1} ^{n-1} \limits
    2^{ d_{n,j} } \cdot
    (d_{n,j})! \cdot 
    (b_{n,j}-d_{n,j})!
    + 2 \cdot
    \prod _{j=1} ^{n-1} \limits
    2^{d_{n-1,j} } \cdot
    (d_{n-1,j})! \cdot
    (b_{n-1,j}-d_{n-1,j})! \\
    & + \sum _{i=3} ^{n-2} \limits
    \prod _{j=1} ^{n-1} \limits
    2^{ d_{i,j}+d_{n+1-i,j} } \cdot
    (d_{i,j}+d_{n+1-i,j})! \cdot 
    (b_{i,j}+b_{n+1-i,j}-d_{i,j}-d_{n+1-i,j})! \\
     = & \ 2 \cdot
     \sum _{i=n-1} ^{n} \limits
    \prod _{j=1} ^{n-1} \limits
    2^{ d_{i,j} } \cdot
    (d_{i,j})! \cdot 
    (b_{i,j}-d_{i,j})! \\
    & + \sum _{i=3} ^{n-2} \limits
    \prod _{j=1} ^{n-1} \limits
    2^{ d_{i,j}+d_{n+1-i,j} } \cdot
    (d_{i,j}+d_{n+1-i,j})! \cdot 
    (b_{i,j}+b_{n+1-i,j}-d_{i,j}-d_{n+1-i,j})!
\end{align*}

\hfill $\square$

\section{Erweiterung der ursprünglichen Höflichkeitsregel}

\noindent Die in Abschnitt 1 angeführte Höflichkeitsregel kann man noch um den Aspekt erweitern, dass die nächste Person immer einen der Plätze auswählt, bei denen möglichst wenige Personen auf den benachbarten Plätzen sitzen:\\

\noindent \textbf{Höflichkeitsregel:} Wer dazu kommt, wählt einen der Plätze so aus, dass der Abstand zu einem besetzten Platz maximal wird. Dabei sollen die Personen die Plätze so besetzen, dass möglichst wenige Personen auf den benachbarten Plätzen sitzen. \\

\noindent Sei $A_n$ die Anzahl der Möglichkeiten, dass $n$ Personen nacheinander einen der $n$ Plätze besetzen und sich dabei an die neue Höflichkeitsregel halten. Damit die neu hinzugefügte Regel relevant wird, muss mindestens eine Situation der Form $OXOX$ vorliegen. Für $n < 4$ kann somit die Formel aus Abschnitt 5 verwendet werden. \\

\noindent \textbf{Behauptung:} Für $n \geq 4$ gilt
\begin{align*}
    A_n &=
    2\cdot 2^{ d_{n,1} } \cdot (d_{n,1})!
    \cdot (b_{n,1}-d_{n,1})! \cdot
    \prod _{j=2} ^{n-1} \limits
    2^{d_{n,j}} \cdot (b_{n,j})! \\
    & + 2 \cdot 2^{d_{n-1,1}} \cdot (d_{n-1,1}+1)!
    \cdot (b_{n-1,1}-d_{n-1,1})! \cdot
    \prod _{j=2} ^{n-1} \limits
    2^{d_{n-1,j}} \cdot (b_{n-1,j})! \\
    & + \sum _{i=3} ^{n-2} \limits
    2^{ d_{i,1}+d_{n+1-i,1} } \cdot (d_{i,1}+d_{n+1-i,1})!
    \cdot (b_{i,1}+b_{n+1-i,1}-d_{i,1}-d_{n+1-i,1})! \ \cdot \\
    &\prod _{j=2} ^{n-1} \limits
    2^{d_{i,j}+d_{n+1-i,j}} \cdot (b_{i,j}+b_{n+1-i,j})!
\end{align*}

\noindent \textbf{Beweis:} Sei $n \geq 4$. Die im Vergleich zu Abschnitt 1 zusätzliche Regel spielt nur beim Besetzen von Plätzen mit Abstand $1$ eine Rolle. \\

\noindent Sei zunächst $i=2$ oder $i=n-1$. Beim Besetzen von Plätzen mit Abstand $1$ werden nun erst von jeder Menge ein Platz und der Platz links oder rechts außen besetzt , da dort nur einer der zwei benachbarten Plätze besetzt ist. Da beide Plätze einer Menge dafür in Frage kommen und es $d_{2,1}+d_{n+1-2,1}$ Mengen gibt, gibt es dafür also $2^{ d_{2,1}+d_{n-1,1}} \cdot (d_{2,1}+d_{n-1,1}+1)!$ Möglichkeiten. Anschließend werden die verbliebenen Plätze mit Abstand $1$ ($b_{2,1}+b_{n-1,1}-d_{2,1}-d_{n-1,1}-1$ als Anzahl der verbliebenen Plätze mit Abstand $1$) besetzt, wofür es $(b_{2,1}+b_{n-1,1}-d_{2,1}-d_{n-1,1}-1)!$ Möglichkeiten gibt.\\

\noindent Bei $i=2$ oder $i=n-1$ gibt es also $2^{ d_{2,1}+d_{n-1,1}} \cdot (d_{2,1}+d_{n-1,1}+1)! \cdot (b_{2,1}+b_{n-1,1}-d_{2,1}-d_{n-1,1}-1)!$ Möglichkeiten, die Plätze mit Abstand $1$ zu besetzen. Für alle anderen Werte von $i$ gibt es nach ähnlicher Herleitung $2^{ d_{i,1}+d_{n+1-i,1} } \cdot (d_{i,1}+d_{n+1-i,1})! \cdot (b_{i,1}+b_{n+1-i,1}-d_{i,1}-d_{n+1-i,1})!$ Möglichkeiten, die Plätze mit Abstand $1$ zu besetzen.\\
    
\noindent Mit den gleichbleibenden Möglichkeiten bei Abständen $ \geq 2$ (siehe dazu Abschnitt 5) folgt dann für $n \geq 4$:
\begin{align*}
    A_n &=
    2\cdot 2^{ d_{1,1}+d_{n,1} } \cdot (d_{1,1}+d_{n,1})!
    \cdot (b_{1,1}+b_{n,1}-d_{1,1}-d_{n,1})! \cdot
    \prod _{j=2} ^{n-1} \limits
    2^{d_{1,j}+d_{n,j}} \cdot (b_{1,j}+b_{n,j})! \\
    & + 2 \cdot 2^{ d_{2,1}+d_{n-1,1} } \cdot (d_{2,1}+d_{n-1,1}+1)!
    \cdot (b_{2,1}+b_{n-1,1}-d_{2,1}-d_{n-1,1}-1)! \ \cdot \\
    &\prod _{j=2} ^{n-1} \limits
    2^{d_{2,j}+d_{n-1,j}} \cdot (b_{2,j}+b_{n-1,j})! \\
    & + \sum _{i=3} ^{n-2} \limits
    2^{ d_{i,1}+d_{n+1-i,1} } \cdot (d_{i,1}+d_{n+1-i,1})!
    \cdot (b_{i,1}+b_{n+1-i,1}-d_{i,1}-d_{n+1-i,1})! \ \cdot \\
    &\prod _{j=2} ^{n-1} \limits
    2^{d_{i,j}+d_{n+1-i,j}} \cdot (b_{i,j}+b_{n+1-i,j})!
\end{align*}

\noindent Mit $d_{1,j}=0, d_{2,j}=0$ und $b_{1,j}=0$ für alle $j$, $b_{2,j}=0$ für $j \geq 2$ und $b_{2,1}=1$ (siehe Behauptungen 3.1, 3.2, 4.1 und 4.2) gilt weiter:
\begin{align*}
    A_n &=
    2\cdot 2^{ d_{n,1} } \cdot (d_{n,1})!
    \cdot (b_{n,1}-d_{n,1})! \cdot
    \prod _{j=2} ^{n-1} \limits
    2^{d_{n,j}} \cdot (b_{n,j})! \\
    & + 2 \cdot 2^{d_{n-1,1}} \cdot (d_{n-1,1}+1)!
    \cdot (b_{n-1,1}-d_{n-1,1})! \cdot
    \prod _{j=2} ^{n-1} \limits
    2^{d_{n-1,j}} \cdot (b_{n-1,j})! \\
    & + \sum _{i=3} ^{n-2} \limits
    2^{ d_{i,1}+d_{n+1-i,1} } \cdot (d_{i,1}+d_{n+1-i,1})!
    \cdot (b_{i,1}+b_{n+1-i,1}-d_{i,1}-d_{n+1-i,1})! \ \cdot \\
    &\prod _{j=2} ^{n-1} \limits
    2^{d_{i,j}+d_{n+1-i,j}} \cdot (b_{i,j}+b_{n+1-i,j})!
\end{align*}

\hfill $\square$

\section{Zusammenfassung}

\noindent Anfänglich wurde das Problem mit $n$ Plätzen und $n$ Personen beschrieben. Besonderes Interesse galt dabei der Anzahl der Möglichkeiten $a_n$. Auch maßgebliche Ideen oder Ansätze wurden in Abschnitt 1 genannt und erläutert. Insbesondere die Anzahl der Personen, die einen Platz mit Abstand $k$ besetzen, und die Anzahl der Mengen waren von entscheidender Bedeutung. Für diese Anzahlen wurden dann später in den Abschnitten 3 und 4 $b_{p,k}$ und $d_{p,k}$ eingeführt. Des Weiteren wurde zwei Schemata für das Hinzufügen von leeren oder besetzten Plätzen in Abschnitt 2 erarbeitet, die in verschiedenen Beweisen in den Abschnitten 3 und 4 verwendet wurden. \\

\noindent Im Abschnitt 3.1 wurde für $k \geq 2$ Folgendes hergeleitet ($m = \lfloor \log_{2}(\frac{p-1}{2k}) \rfloor$ für $p \geq 5$):

\begin{equation*}
    b_{p,k} = 
    \begin{cases}
        0
        \ \ & \text{für} \ \ p < k+1  \\
        
        1
        \ \ & \text{für} \ \ p = k+1  \\
        
        0
        \ \ & \text{für} \ \ k+1<p < 1+2k  \\
    
        2^m 
        \ \ & \text{für} \ \ 1+2^{m}\cdot 2k \leq p \leq 1+2^{m}(2k+1) \\

        1+2^m(2k+2)-p
        \ \ & \text{für} \ \ 1+2^{m}(2k+1) < p \leq 1+2^{m}(2k+2)\\
        
        0  
        \ \ & \text{für} \ \ 1+2^{m}(2k+2)< p \leq 1+2^{m}(4k-2) \\
        
        p-1-2^{m}(4k-2)
        \ \ & \text{für} \ \ 1+2^{m}(4k-2) < p < 1+2^{m+1}2k

    \end{cases}
\end{equation*}

\noindent Eine Möglichkeit, $b_{p,1}$ zu berechnen, wurde dann in Abschnitt 3.2 gefunden ($m = \lfloor \log_{2}(\frac{p-1}{3}) \rfloor$ für $p \geq 4$):

\begin{equation*}
    b_{p,1} =
    \begin{cases}
        0 \ \ & \text{für} \ \ p=1 \\
        1 \ \ & \text{für} \ \ p=2 \ \text{oder} \ p=3 \\
        2^{m+1} \ \ & \text{für} \ \ 1+2^{m} \cdot 3 \leq p \leq 1+2^{m} \cdot 4 \\
        p-1-2^{m+1} \ \ & \text{für} \ \ 1+2^{m} \cdot 4 < p < 1+2^{m+1}\cdot 3
    \end{cases}
\end{equation*}

\noindent Weiter kann $d_{p,k}$ nach Abschnitt 4.1 für $k \geq 2$ folgendermaßen berechnet werden ($m = \lfloor \log_{2}(\frac{p-1}{2k}) \rfloor$ für $p \geq 5$):

\begin{equation*}
    d_{p,k} =
    \begin{cases}
        0 
        & \text{für} \ \  p < 1+2k \\

        p - 1-2^{m}\cdot 2k
        & \text{für} \ \ 
        1+2^{m}\cdot 2k \leq p \leq 1+2^{m}(2k+1) \\

        1+2^{m}(2k+2)-p
        & \text{für} \ \
        1+2^{m}(2k+1) < p \leq 1+2^{m}(2k+2) \\

        0 
        & \text{für} \ \
        1+2^{m}(2k+2) < p < 1+2^{m+1} \cdot 2k
    \end{cases}
\end{equation*}

\noindent Mit diesen Möglichkeiten, $b_{p,k}$ und $d_{p,k}$ zu berechnen, den anfänglichen Ideen und $d_{p,1}=0$ für alle $p$ ist dann in Abschnitt 5 eine Formel für $a_n$ hergeleitet worden:

\begin{equation*}
    a_n = \sum ^{n} _{i=1} \limits \prod \limits ^{n-1} _{j=1} 2^{d_{i,j} + d_{n+1-i,j}} \cdot (b_{i,j} + b_{n+1-i,j})!
\end{equation*}

\noindent Des Weiteren wurde $a_n$ in Abschnitt 6 für $n \geq 2$ nach unten ($m_i = \left \lfloor\log_2 (\frac{i-1}{4}) \right \rfloor$) und für alle $n$ nach oben abgeschätzt:
\begin{align*}
     & 2 \cdot \sum 
    ^{n} 
    _{i=\left \lceil \frac{n}{2} \right \rceil +1} \limits 
    \left (\left \lfloor \frac{i-1}{2} \right \rfloor
     + \left \lfloor \frac{n-i}{2} \right \rfloor \right )! 
    \cdot
    \prod \limits ^{m_i} _{j=1} (2^{j-1})!
    \leq
    a_n \\
    & a_n 
    \leq 
    \sum _{i=1} ^n \limits
    \left ( \left \lceil \frac{2(i-1)}{3} \right \rceil +\left \lceil \frac{2(n-i)}{3} \right \rceil \right )! \cdot
    \left ( n-\left \lceil \frac{2(i-1)}{3} \right \rceil-\left \lceil \frac{2(n-i)}{3} \right \rceil \right )!
\end{align*}

\noindent Anschließend wurde dann in Abschnitt 7 die Formel für $a_n$ mit den Abschätzungen verglichen. In Abschnitt 8 konnten Formeln für die OEIS-Folgen A166079, A095236, A095240 und A095912 angegeben und bewiesen werden. Schlussendlich wurde die Höflichkeitsregel aus Abschnitt 1 in Abschnitt 9 erweitert. Auch dort konnte eine Formel für $A_n$ als Anzahl der Möglichkeiten hergeleitet werden.

\newpage

\section{Literaturverzeichnis}

[1] \ Johann Beurich alias DorFuchs: \textit{Das Pissoir-Problem}, 2022, \\
URL: \url{https://www.youtube.com/watch?v=a36zHPlbd2g} \\

\noindent [2] The On-Line Encyclopedia of Integer Sequences (OEIS), \url{https://oeis.org}, Folge A166079 \\

\noindent [3] The On-Line Encyclopedia of Integer Sequences (OEIS), \url{https://oeis.org}, Folge A095236 \\

\noindent [4] The On-Line Encyclopedia of Integer Sequences (OEIS), \url{https://oeis.org}, Folge A095240 \\

\noindent [5] The On-Line Encyclopedia of Integer Sequences (OEIS), \url{https://oeis.org}, Folge A095912\\

\vspace{1cm}

\noindent Weitere Literatur im Kontext der Problemstellung:
\begin{itemize}
    \item  Evangelos Kranakis und Danny Krizanc: \textit{The Urinal Problem}, 2010, \\ URL: \url{https://people.scs.carleton.ca/~kranakis/Papers/urinal.pdf}
    \item The On-Line Encyclopedia of Integer Sequences (OEIS), \url{https://oeis.org}, Folge A358056
\end{itemize}

\end{document}